\documentclass [10pt,a4paper]{article}
\usepackage{amssymb}
\usepackage{enumerate,verbatim}
\usepackage{amsmath}
\usepackage{amsfonts,mathrsfs}
\usepackage{amsthm,wasysym}
\usepackage{epsfig,lpic,tikz}
\usepackage[applemac]{inputenc}
\usepackage[english]{babel}

\usepackage[width=0.9\textwidth]{caption}

\usepackage{pdfsync}
\usepackage{xcolor,hyperref}
\usepackage[T1]{fontenc}
\usepackage{lpic}
\usepackage{mathtools}

\hypersetup{
    colorlinks=true,
    linkcolor=blue,
    filecolor=magenta,      
    urlcolor=blue,
    pdftitle={Overleaf Example},
    pdfpagemode=FullScreen,
    }

\definecolor{green}{rgb}{0,0.5,0}
\definecolor{teal}{RGB}{38,140,83}

\usepackage{marvosym}

\newcounter{teoremaganso}
\newcounter{appendix}

\newcounter{coryganso}

\renewenvironment{abstract}{\small\quotation\noindent
 {\bfseries \abstractname .}}{\endquotation \par}

\newenvironment{trivlistparm}[2]{\trivlist
\item[{#1 #2}]}{\endtrivlist}

\newenvironment{prooftext}[1]{\trivlistparm{\bfseries}{#1}}{\Qed\endtrivlistparm}
\newenvironment{prova}{\trivlistparm{\bfseries}{Proof.}}{\Qed\endtrivlistparm}

\catcode`\@=11

\def\resetthefootnote{\renewcommand{\thefootnote}{\@arabic\c@footnote} }
\def\@principiremex#1{\trivlist
 \item[\hskip \labelsep{\bfseries #1\ \thetheo.}]\ignorespaces}
\def\opar@principiremex#1[#2]{\trivlist
 \item[\hskip \labelsep{\bfseries #1\ \thetheo\ (#2).}]\ignorespaces}

\newcommand{\newTHEOremrom}[2]{\newenvironment{#1}{\refstepcounter{theo}\@ifnextchar[{\opar@principiremex{#2}}
{\@principiremex{#2}}}{\qedB\endtrivlist}} \catcode`\@=12
\DeclareMathSymbol{\square}{\mathord}{AMSa}{"03}
\newcommand{\qedB}{\nopagebreak\hspace*{\fill}$\square$\par}
\newcommand{\Qed}{\nopagebreak\hspace*{\fill}{\vrule width6pt height6pt depth0pt}\par}

\newtheorem {theo} {Theorem} [section]

\newtheorem {lem} [theo] {Lemma}
\newtheorem {bigtheo} [teoremaganso] {Theorem}

\newTHEOremrom {defi} {Definition}
\newTHEOremrom {obs} {Remark}
\newTHEOremrom {ex} {Example}

\newcommand{\refc}[1]{\mbox{$(\ref{#1})$}}

\newcommand{\teoc}[1]{Theorem~\ref{#1}}

\usepackage{amsmath,amssymb}
\usepackage{graphicx}
\usepackage{color}
\usepackage{psfrag} 
\usepackage{overpic} 
\usepackage{array} 

\title{\bf The period of the limit cycle bifurcating \\ from a persistent  polycycle
}

\author{D. Mar\'{\i}n, L. Queiroz and J. Villadelprat
\\*[.1truecm]
{\small \textsl{Departament de Matem{\`a}tiques, Edifici Cc,
Universitat Aut{\`o}noma de Barcelona,}}\\*[-.05truecm]
{\small\textsl{08193 Cerdanyola del Vall\`es (Barcelona), Spain}}
\\*[-.05truecm]
{\small \textsl{Centre de Recerca Matem\`atica, Edifici Cc, Campus de Bellaterra,}}\\*[-.05truecm]
{\small \textsl{08193 Cerdanyola del Vall\`es (Barcelona), Spain}}
\\*[.1truecm]
{\small \textsl{Universidade Estadual Paulista (UNESP), Instituto de Bioci\^encias, Letras e Ci\^encias Exatas,}}
\\*[-.05truecm]
{\small \textsl{R. Cristov\~ao Colombo, 2265, 15054-000 S\~ao Jos\'e do Rio Preto, SP, Brazil}}
\\*[.1truecm]
{\small \textsl{Departament d'Enginyeria Inform{\`a}tica i Matem{\`a}tiques, ETSE,}}
\\*[-.05truecm]
{\small \textsl{Universitat Rovira i Virgili, 43007 Tarragona, Spain}}}
                  
\date{\today}

\begin{document}

\maketitle

\begin{abstract}
We consider smooth families of planar polynomial vector fields $\{X_\mu\}_{\mu\in\Lambda}$, where $\Lambda$ is an open subset of $\mathbb{R}^N$, for which there is a hyperbolic polycycle $\Gamma$ that is persistent (i.e., such that none of the separatrix connections is broken along the family). It is well known that in this case the cyclicity of~$\Gamma$ at~$\mu_0$ is zero unless its graphic number $r(\mu_0)$ is equal to one. It is also well known that if $r(\mu_0)=1$ (and some generic conditions on the return map are verified) then the cyclicity of $\Gamma$ at $\mu_0$ is one, i.e., exactly one limit cycle bifurcates from $\Gamma$. In this paper we prove that this limit cycle approaches $\Gamma$ exponentially fast and that its period goes to infinity as $1/|r(\mu)-1|$ when $\mu\to\mu_0.$ Moreover, we prove that if those generic conditions are not satisfied, although the cyclicity may be exactly 1, the behavior of the period of the limit cycle is not determined.
\end{abstract}

\section{Introduction and main results}

This work deals with the study of the period of limit cycles arising in bifurcations for families of smooth planar vector fields $\{X_\mu\}_{\mu\in\Lambda}$, where $\Lambda\subset\mathbb{R}^N$. From the classification of first-order structurally unstable vector fields (see, for instance \cite{Andronov, Kocak, Soto74}), the generic compact isolated bifurcations for one-parameter (i.e., $N=1$) families of smooth planar vector fields which give rise to periodic orbits for $\mu\to\mu_0$ are: the Andronov-Hopf bifurcation, the bifurcation of a semi-stable periodic orbit, the saddle-node loop and the saddle loop bifurcations. These are referred to as the \emph{elementary bifurcations}. In \cite{GasManVil05} the authors determined the behavior of the period $\mathscr{T}(\mu)$ of the limit cycle of $X_\mu$ arising from a elementary bifurcation as $\mu\to \mu_0$. More precisely, they obtained the principal term of the expression of $\mathscr{T}(\mu)$ which is comprised in the following list:
\begin{enumerate}[$(i)$]
	\item $\mathscr{T}(\mu)\sim T_0+T_1(\mu-\mu_0)$ for the Andronov-Hopf bifurcation;
	\item  $\mathscr{T}(\mu)\sim T_0+T_1\sqrt{|\mu-\mu_0|}$ for the bifurcation from a semi-stable periodic orbit;
	\item $\mathscr{T}(\mu)\sim T_0/\sqrt{|\mu-\mu_0|}$ for the saddle-node loop bifurcation;
	\item $\mathscr{T}(\mu)\sim T_0\log|\mu-\mu_0|$ for the saddle loop bifurcation.
\end{enumerate}
Here $T_0$ and $T_1$ are constants with $T_0\neq 0$ and we use the notation $\mathscr{T}(\mu)\sim a+f(\mu-\mu_0)$ as $\mu$ tends to $\mu_0$ meaning that $\lim_{\mu\to \mu_0}(\mathscr{T}(\mu)-a)/f(\mu-\mu_0)=1$. Accordingly the principal term of the period of the periodic orbit arising from generic elementary bifurcations characterizes the bifurcation. In what concerns to applications, the information concerning $\mathscr{T}(\mu)$ can be useful to estimate parameters associated to a system, for instance, when studying neuron activities in the brain with the aim of determining the synaptic conductance that it receives (see \cite{Toni}).

When studying the case $(iv)$ above, the authors in \cite{GasManVil05} assume the breaking of the saddle loop connection, which is ``the first order condition'' for the bifurcation of a limit cycle to occur.
The present work started as a natural continuation of this previous study by considering the case in which the homoclinic connection remains unbroken and the parameter space is not necessarily one-dimensional. 
That being said, the tools and techniques applied in our approach enabled us to extend the results obtained for the saddle loop case to more general situations, namely to any hyperbolic polycycle that is persistent (i.e. such that none of the separatrix connections is broken along the perturbation). For a persistent polycycle $\Gamma$, we have an associated return map defined from a fixed transversal section $\Sigma$ parametrized smoothly by $\sigma(s;\mu)$ which writes as
\begin{equation}\label{eqRintro}
\mathscr{R}(s;\mu)=s^{r(\mu)}\left(A(\mu)+R(s;\mu)\right),
\end{equation}
where $r(\mu)$ is the product of the hyperbolicity ratios of all saddles in $\Gamma$ (also called the \emph{graphic number}~\cite{GasMan2002}) and $R(s;\mu)$ is a well-behaved remainder. It is well known (see \cite{Soto81}) that the stability of $\Gamma$ is determined by the graphic number in case that $r(\mu_0)\neq 1$, which is the case considered in \cite{GasManVil05}. There are however other situations that may occur. The list, from the most to the least generic case, is the following:
\begin{enumerate}[$(a)$]
\item A saddle connection breaks and $r(\mu_0)\neq 1$;
\item A saddle connection breaks and $r(\mu_0)=1$;
\item  No saddle connection breaks and $r(\mu_0)\neq 1$;
\item  No saddle connection breaks and $r(\mu_0)=1$.
\end{enumerate}
As we  already said, case $(a)$ is treated in \cite{GasManVil05}, whereas there is no bifurcation of limit cycle in case $(c)$, see Remark~\ref{obsrneq1}.
The present paper is addressed to case $(d)$ and we have obtained two main results. The first one is the following.

\begin{bigtheo}\label{TeoPeriodmPoly}
Let $\{X_\mu\}_{\mu\in\Lambda}$ be a smooth family of planar polynomial vector fields. Let us fix $\mu_0\in\Lambda$ such that $X_{\mu_0}$ has a persistent hyperbolic polycycle $\Gamma$.
If $\Gamma\cap\ell_\infty\neq\emptyset$ we assume additionally that
 the infinite line~$\ell_\infty$ is invariant for all $\mu\in\Lambda$. Fix a transversal section $\Sigma$ to $\Gamma$ parametrized smoothly by $s\mapsto\sigma(s;\mu)$ with $\sigma(0;\mu)\in\Gamma$ for all $\mu\in\Lambda$. Then the following holds:
	\begin{itemize}
		\item[$(a)$] The associated return map writes as
		$$\mathscr{R}(s;\mu)=s^{r(\mu)}\left(A(\mu)+\mathcal{F}^{\infty}_{\epsilon}(\mu_0)\right),$$
		for any $\epsilon>0$ small enough, where $r(\mu)$ is the product of all the hyperbolicity ratios of the saddles in~$\Gamma$ and $A\in \mathscr C^\infty(\Lambda)$ with $A(\mu_0)>0$;
		\item[$(b)$] The associated return time is given by
		$$\mathcal{T}(s;\mu)=\bar{T}_0(\mu)\log s +\mathcal{F}^{\infty}_{0}(\mu_0),$$
		where $\bar{T}_{0}\in\mathscr{C}^{\infty}(\Lambda)$. If at least one of the saddles in $\Gamma$ is finite then $\bar{T}_{0}(\mu)<0$.
	\end{itemize}
Moreover, if $r(\mu)-1$ changes sign at $\mu_0$ and $A(\mu_0)\neq 1$ then the following also holds:
	\begin{itemize}
		\item[$(c)$] $\rm{Cycl}(\Gamma,\mu_0)=1$ and the limit cycle $\gamma_\mu$ that bifurcates the polycycle approaches it at least exponentially fast, i.e., there exists $K>0$ such that 
		\[
		\lim_{\mu\to\mu_0}\frac{s(\mu)}{e^{\frac{-K}{|1-r(\mu)|}}}=0,
		\] 
		where $s(\mu)$ is the location of $\gamma_\mu$ at $\Sigma$, i.e., the unique positive solution of $\mathscr{R}(s;\mu)=s$;
		\item[$(d)$] If at least one of the saddles at $\Gamma$ is finite then the period $\mathscr{T}(\mu)$ of $\gamma_\mu$ goes to infinity as $\frac{1}{|1-r(\mu)|}$ i.e. $\lim\limits_{\mu\to\mu_0} |1-r(\mu)|\mathscr{T}(\mu)$ is a real positive number.
	\end{itemize}
\end{bigtheo} 

The notation $\mathcal{F}_{L}^{K}(\mu_0)$ in the remainder refers to the notion of \emph{finitely flat functions} (see Definition~\ref{defiFlat}), which is the new information we add to the known results given in $(a)$ and $(b)$. We include them for completeness and to state precisely the hypothesis in $(c)$ and $(d)$, which to the best of our knowledge are new.

Under the hypothesis of Theorem \ref{TeoPeriodmPoly}, if the parameter space is one-dimensional (i.e., $N=1$) then the graphic number writes as $r(\mu)=1+\tfrac{1}{T_0}(\mu-\mu_0)^k+O(\mu^{k+1})$ for an odd positive integer $k$. Thus, by applying Theorem \ref{TeoPeriodmPoly}, the behavior of the period of the limit cycle bifurcating from a persistent polycycle is $\mathscr{T}(\mu)\sim T_0/|\mu-\mu_0|^k$, which is an entirely different behavior from any of the elementary bifurcations discussed above.
However, it coincides with the behavior of the generalized Andronov-Hopf bifurcation considered in
 \cite[Theorem~7]{GasManVil05}. Thus, one cannot classify bifurcations of limit cycles based uniquely on the behavior of their periods.

In the case which the above conditions on $r(\mu)$ and $A(\mu)$ are not satisfied, although we may have cyclicity exactly 1, the behavior of the period of the bifurcating limit cycle is not determined. In fact, we solve a sort of inverse problem when $r(\mu_0)=A(\mu_0)=1$ and the gradients of $A$ and $r$ are linearly independent at~$\mu_0$. In this regard, our second main result is the following.

\begin{bigtheo}\label{TeoCenter}
Let $\{X_\mu\}_{\mu\in\Lambda}$ be a smooth family of planar polynomial vector fields. Let us fix $\mu_0\in\Lambda$ such that $X_{\mu_0}$ has a persistent hyperbolic polycycle $\Gamma$ with at least a finite saddle. If $\Gamma\cap\ell_\infty\neq\emptyset$ we assume additionally that 
 the infinite line $\ell_\infty$ is invariant for all $\mu\in\Lambda$. Fix a transversal section $\Sigma$ to $\Gamma$ parametrized smoothly by $s\mapsto\sigma(s;\mu)$ with $\sigma(0;\mu)\in\Gamma$ for all $\mu\in\Lambda$. Then the associated return map writes as
	$$\mathscr{R}(s;\mu)=s^{r(\mu)}\left(A(\mu)+\mathcal{F}^{\infty}_{\epsilon}(\mu_0)\right),$$
	for any $\epsilon>0$ small enough, where $r(\mu)$ is the product of all the hyperbolicity ratios of the saddles in $\Gamma$ and $A\in \mathscr C^\infty(\Lambda)$ with $A(\mu_0)>0$. Suppose that $A(\mu_0)=r(\mu_0)=1$ and ${\rm rank}\{\nabla A,\nabla r\}\vert_{\mu=\mu_0}=2$.
	Then for every function $\tau\in\mathscr{C}^{1}(0,\delta)$ satisfying $\tau(\alpha)>0$ for $\alpha>0$ and 
\begin{equation}\label{cond}
\lim\limits_{\alpha\to 0^{+}}\alpha\tau(\alpha)=\lim\limits_{\alpha\to 0^{+}}\alpha^2 \tau'(\alpha)=\lim\limits_{\alpha\to 0^{+}}\frac{\log\alpha}{\tau(\alpha)}=0,
\end{equation}
there exists a differentiable arc at $\mu_0$ in the parameter space $\Lambda$, i.e., a $\mathscr C^1$ map $m:(0,\delta_0)\to\Lambda$ with $\lim_{\alpha\to 0^{+}}m(\alpha)=\mu_0$, such that $X_{m(\alpha)}$ has a limit cycle that approaches the polycycle $\Gamma$ as $\alpha\to 0^+$ with period $\mathscr{T}(\alpha)$ verifying that $\lim\limits_{\alpha\to 0^+}\mathscr{T}(\alpha)/\tau(\alpha)$ is a positive real number.
\end{bigtheo} 
For instance, the function $\tau(\alpha)=\alpha^l$  fulfills conditions (\ref{cond}) provided that $l\in(-1,0)$.

This paper is structured as follows. In Section \ref{SecDefi}, we present the preliminary definitions necessary for the development of our work, including the definition of the finitely flat functions and their properties. Additionally, in this section, we prove the technical results that we use in the localization of the limit cycles bifurcating from hyperbolic persistent polycycles. Section \ref{SecTeoA} is dedicated to the proof of Theorem \ref{TeoPeriodmPoly}, which deals with the period of the unique limit cycle bifurcating from a persistent polycycle in the case $r(\mu_0)=1$ and $A(\mu_0)\neq 1$. In Section \ref{SecTeoB}, we prove Theorem \ref{TeoCenter}, which addresses the case $r(\mu_0)=A(\mu_0)=1$ in which the behavior of the period of limit cycles is undetermined. In the Appendix we gather some previous results about the Dulac map and time which are used throughout the present text.

\section{Definitions and preliminary results}\label{SecDefi}

We now define precisely the notions of persistent polycycle and the cyclicity of a polycycle.

	\begin{defi}
		Let $X$ be a two-dimensional vector field. A \emph{graphic} $\Gamma$ for $X$ is a compact non-empty invariant subset which is a continuous image of $\mathbb{S}^1$ and consists of a finite number of (not necessarily distinct) isolated singular points $\{p_1,\dots,p_m,p_{m+1}=p_1\}$ and compatibly oriented separatrices $\{s_1,\dots,s_m\}$ connecting them (meaning that $s_i$ has $\{p_i\}$ as the $\alpha$-limit set and $\{p_{i+1}\}$ as the $\omega$-limit set). A graphic for which all its singular points are hyperbolic saddles is said to be \emph{hyperbolic}. A \emph{polycycle} is a graphic with a return map defined on one of its sides.
	\end{defi}

\begin{defi}[Persistent polycycle]
Let $\{X_\mu\}_{\mu\in\Lambda}$ be a smooth family of planar vector fields such that $\Gamma$ is a hyperbolic polycycle of $X_{\mu_0}$. We say that $\Gamma$ is a \emph{persistent polycycle} when all of its separatrix connections remain unbroken inside the family $\{X_\mu\}_{\mu\in\Lambda}$.
	\end{defi}
	
	\begin{defi}
		Let $\{X_\mu\}_{\mu\in\Lambda}$ be a family of vector fields on $\mathbb{S}^2$ and suppose that $\Gamma$ is a polycycle for $X_{\mu_0}$. We say that $\Gamma$ has finite cyclicity in the family $\{X_\mu\}_{\mu\in\Lambda}$ if there exist $\kappa\in\mathbb{N}$, $\varepsilon>0$ and $\delta>0$ such that any $X_{\mu}$ with $\vert\mu-\mu_0\vert<\delta$ has at most $\kappa$ limit cycles $\gamma_i$ with $\rm{dist}_H(\Gamma,\gamma_i)<\varepsilon$ for $i=1,\dots,\kappa$. The minimum of such $\kappa$ when $\delta$ and $\varepsilon$ go to zero is called \emph{cyclicity} of $\Gamma$ in  $\{X_\mu\}_{\mu\in\Lambda}$ at $\mu=\mu_0$ and denoted by $\rm{Cycl}(\Gamma,\mu_0)$.
\end{defi}

We introduce the notion of finitely flat functions that play a substantial role in the results of this paper.

\begin{defi}
	Consider $K\in\mathbb{Z}_{\geqslant 0}\cup\{\infty\}$ and an open set $U\subset\mathbb{R}^{N}$. We say that a function $\psi(s;\mu)$ belongs to class $\mathscr C^K_{s>0}(U)$, if there exists an open neighborhood $\Omega$ of ${0}\times U$ in $\mathbb{R}^{N+1}$ such that $(s;\mu)\mapsto \psi(s,\mu)$ is $\mathscr C^K$ on $\Omega\cap\left((0,+\infty)\times U\right)$.
\end{defi}

\begin{defi}[Finitely flat functions]\label{defiFlat}
	Consider $K\in\mathbb{Z}_{\geqslant 0}\cup\{\infty\}$ and an open set $U\subset\mathbb{R}^{N}$. Given $L\in\mathbb{R}$ and $\mu_0\in U$, we say that $\psi(s;\mu)\in \mathscr C^K_{s>0}(U)$ is $(L,K)$-flat with respect to $s$ at $\mu_0$, and we write $\psi\in\mathcal{F}^{K}_{L}(\mu_0)$, if for each $\nu=(\nu_0,\dots,\nu_{N})\in\mathbb{Z}^{N+1}_{\geqslant 0}$ with $|\nu|\leqslant K$, there exist a neighborhood $V$ of $\mu_0$ and $C$, $s_0>0$ such that
	$$\partial_\nu\psi(s;\mu):=\left\vert\dfrac{\partial^{|\nu|}\psi(s;\mu)}{\partial s^{\nu_0}\partial\mu_1^{\nu_1}\cdots\partial\mu_N^{\nu_N}}\right\vert\leqslant Cs^{L-\nu_0}\;\text{for all } s\in(0,s_0)\text{ and } \mu\in V.$$
	If $W$ is a (not necessarily open) subset of $U$, then $\mathcal{F}^{\infty}_{L}(W)=\bigcap\limits_{\mu_0\in W}\mathcal{F}^{\infty}_{L}(\mu_0)$. 
\end{defi}

The next two results were proven in \cite{MarVilDulacLocal} and will be very useful to work with finitely flat functions.

\begin{lem}\label{LemaPropFlat}
	Let $U$ and $U'$ be open sets of $\mathbb{R}^N$ and $\mathbb{R}^{N'}$ respectively and consider $W\subset U$ and $W'\subset U'$. Then, the following holds:
	\begin{itemize}
		\item[(a)] $\mathcal{F}^{K}_{L}(W)\subset\mathcal{F}^{K}_{L}(\hat{W})$ for any $\hat{W}\subset W$;
		\item[(b)]$\mathcal{F}^{K}_{L}(W)\subset\mathcal{F}^{K}_{L}(W\times W')$;
		\item[(c)] $\mathscr C^K(U)\subset\mathcal{F}^{K}_{0}(W)$;
		\item[(d)] If $K\geqslant K'$ and $L\geqslant L'$ then $\mathcal{F}^{K}_{L}(W)\subset \mathcal{F}^{K'}_{L'}(W)$;
		\item[(e)] $\mathcal{F}^{K}_{L}(W)$ is closed under addition;		\item[(f)] If $f\in\mathcal{F}^{K}_{L}(W)$ and $\nu\in\mathbb{Z}^{N+1}_{\geqslant 0}$ with $|\nu|\leqslant K$ then $\partial_{\nu}f\in\mathcal{F}^{K-|\nu|}_{L-\nu_0}(W)$;
		\item[(g)] $\mathcal{F}^{K}_{L}(W)\cdot\mathcal{F}^{K'}_{L'}(W)\subset \mathcal{F}^{K}_{L+L'}(W)$;
		\item [(h)] Assume that $\phi:U'\to U$ is a $\mathscr C^K$ function with $\phi(W')\subset W$ and let us take $g\in \mathcal{F}^{K}_{L'}(W')$ with $L'>0$ and verifying $g(s;\eta)>0$ for all $\eta\in W'$ and $s>0$ small enough. Consider also any $f\in\mathcal{F}^{K}_{L}(W)$. Then $h(s;\eta):=f(g(s;\eta);\phi(\eta))$ is a well-defined function that belongs to $\mathcal{F}^{K}_{LL'}(W')$.
	\end{itemize}
\end{lem}

\begin{lem}\label{LemaExt}
Let $U$ be an open set of $\mathbb{R}^N$, $K\in\mathbb{Z}_{\geqslant 0}$ and $g(s;\mu)\in\mathscr{C}^K_{s>0}(U)$ such that, for some $W\subset U$ and $L\in\mathbb{R}$, $g(s;\mu)\in\mathcal{F}^K_L(W)$. If $L>K$, then $g$ extends to a $\mathscr{C}^K$ function $\tilde{g}$, defined in some open neighborhood of $\{0\}\times W\in\mathbb{R}^{N+1}$, satisfying $\partial_\nu\tilde{g}(0;\mu)=0$ for all $\mu\in W$ and $\nu\in\mathbb{Z}_{\geqslant 0}^{N+1}$ with $|\nu|\leqslant K$.
\end{lem}

In what follows, we present several technical results which will be applied in the study of limit cycles bifurcating from hyperbolic persistent polycycles. These are essentially the machinery behind Theorems~\ref{TeoPeriodmPoly} and \ref{TeoCenter}.
\begin{lem}\label{LemaPoly}
	Let $f_1(\mu),\dots,f_m(\mu)$ be $\mathscr C^{k+1}$ functions on an open subset $\Lambda\subset\mathbb{R}^N$ and $k\in\mathbb{Z}_{\geqslant 0}$. For each $j=1,\ldots,m$,
	let $\partial^{\leqslant k}_\mu f_j(\mu)$ denote the vector whose $\binom{N+k}{k}$ entries are all the partial derivatives of order less or equal to $k$ of  $f_j(\mu)$. Then 
for every $P\in\mathbb{R}[x_1,\dots,x_{m\binom{N+k}{k}}]$ and each $i\in\{1,\dots,N+1\}$ there exists $Q_i\in\mathbb{R}[x_1,\dots,x_{m\binom{N+k+1}{k+1}}]$ such that
if $G(\mu)=P(\partial^{\leqslant k}_\mu f_1,\dots,\partial^{\leqslant k}_\mu f_m)$ then
	$\partial_{\mu_i}G(\mu)=Q_i(\partial^{\leqslant k+1}_\mu f_1,\dots,\partial^{\leqslant k+1}_\mu f_m).$
\end{lem}

The above lemma is a direct consequence of the chain rule and the fact that the ring of polynomials is closed under differentiation.

\begin{lem}\label{LemaSFlat}
	For any $L>0$ and $K\in\mathbb{Z}_{\geqslant 0}$, the function $s(\alpha;\mu,z)=c(\mu)^{-1/\alpha}(1+z)$, with $c(\mu)$ being a  $\mathscr C^\infty$ function an open subset $U\subset\mathbb{R}^N$ and $c(\mu_0)>1$ belongs to $\mathcal{F}_L^{K}(\mu_0,0)$.
\end{lem}
\begin{prova}
	First we compute $\lim\limits_{\alpha\to 0^{+}}\alpha^{-L}c(\mu)^{-1/\alpha}$. Since $c(\mu_0)>1$, there is an open neighborhood of $U'$ containing $\mu_0$ such that $c(\mu)>\delta>1$ for all $\mu\in U'$. Moreover,
	$$\log\left(\dfrac{c(\mu)^{-1/\alpha}}{\alpha^L}\right)=-\frac{1}{\alpha}\log c(\mu)-L\log\alpha=-\frac{1}{\alpha}\left(\log c(\mu)-L\alpha\log\alpha\right),$$
	and then, we have that $\lim\limits_{\alpha\to 0^{+}}\log\left(\alpha^{-L}c(\mu)^{-1/\alpha}\right)=-\infty$ and therefore $\lim\limits_{\alpha\to 0^{+}}\alpha^{-L}c(\mu)^{-1/\alpha}=0$. Since $c(\mu)>\delta>1$, these limits are uniform and hence, there exist $C>0$ and a neighborhood $V$ of $(\mu_0,0)$ such that $|s(\alpha;\mu,z)|\leqslant C\alpha^L$ for $\alpha$ small enough and $(\mu,z)\in V$. Thus, $s(\alpha;\mu,z)\in\mathcal{F}_{L}^{0}(\mu_0,0)$.
	
	Now, to see that $s(\alpha;\mu,z)\in\mathcal{F}_{L}^{1}(\mu_0,0)$, it is sufficient to verify that the first derivatives of $s(\alpha;\mu,z)$ are bounded in a neighborhood of $(\mu_0,0)$. This is a direct consequence of their expressions, given as follows:
	\begin{eqnarray}\label{eqPartialsmu}
		\dfrac{\partial s(\alpha;\mu,z)}{\partial z}&=&c(\mu)^{-1/\alpha},\nonumber\\
		\dfrac{\partial s(\alpha;\mu,z)}{\partial \alpha}&=&c(\mu)^{-1/\alpha}\alpha^{-2}(1+z)\log c(\mu),\\
		\dfrac{\partial s(\alpha;\mu,z)}{\partial \mu_i}&=&-\frac{1}{\alpha}c(\mu)^{-1/\alpha-1}(1+z)\frac{\partial c}{\partial \mu_i}(\mu)=-\frac{1}{\alpha}c(\mu)^{-1/\alpha}(1+z)\frac{\partial c}{\partial \mu_i}(\mu)c(\mu)^{-1}.\nonumber
	\end{eqnarray}
	Thus, we have
	$$\lim\limits_{\alpha\to 0^{+}}\alpha^{-L}\left\vert\dfrac{\partial s(\alpha;\mu,z)}{\partial \alpha^{\nu_0}\partial\mu_1^{\nu_1}\cdots\partial\mu_N^{\nu_N}\partial z^{\nu_{N+1}}}\right\vert=0,$$
	uniformly, for $\nu_0+\dots+\nu_{N+1}=1$. Therefore $s(\alpha;\mu,z)\in\mathcal{F}_{L}^1(\mu_0,0)$.
	\medskip
	
	We then proceed by induction to prove the following claim:
	
	\textbf{Claim.} For any $\nu\in\mathbb{Z}_{\geqslant 0}^{N+2}$ with $|\nu|=k\leqslant K$, we have
	\begin{eqnarray}
		\partial_{\nu}s(\alpha;\mu,z)&=&\dfrac{\partial^{|\nu|} s(\alpha;\mu,z)}{\partial \alpha^{\nu_0}\partial\mu_1^{\nu_1}\cdots\partial\mu_N^{\nu_N}\partial z^{\nu_{N+1}}}=\alpha^{-2k}c(\mu)^{-1/\alpha}P_{\nu}\left(\log c(\mu),\alpha,z,c(\mu)^{-1},\partial^{\leqslant k}_{\mu}c(\mu)\right),\nonumber
	\end{eqnarray}
	where $P_{\nu}$ is a polynomial with real coefficients and $\partial^{\leqslant k}_{\mu}c(\mu)$ denotes a vector with all the partial derivatives of $c(\mu)$ of order less or equal to $k$ on $\mu$. 
	
	For $k=1$, by equation \eqref{eqPartialsmu}, this is true. Now, assume it holds for $k=\tilde{k}\geqslant 1$. For $\nu\in\mathbb{Z}_{\geqslant 0}^{N+2}$ with $|\nu|=\tilde{k}+1$, we must consider one of the following cases:
	
	\textbf{1. }$\nu_{N+1}\geqslant 2$. In this case, by Schwarz Theorem, $\partial_{\nu}s(\alpha;\mu,z)=0$ and the claim holds.
	
	\textbf{2. }$\nu_{0}\geqslant 1$. In this case, $\partial_{\nu}s(\alpha;\mu,z)=\partial_\alpha\partial_{\nu-e_0}s(\alpha;\mu,z)$, where $e_0=(1,0,\dots,0)\in\mathbb{Z}_{\geqslant 0}^{N+2}$. Thus, since $|\nu-e_0|=\tilde{k}$, we have
{\small\begin{eqnarray}
			\partial_\alpha\partial_{\nu-e_0}s(\alpha;\mu,z)&=&\partial_\alpha\left(\alpha^{-2\tilde{k}}c(\mu)^{-1/\alpha}P_{\nu-e_0}\left(\log c(\mu),\alpha,z,c^{-1}(\mu),\partial^{\leqslant \tilde{k}}_{\mu}c(\mu)\right)\right)=\nonumber\\
			&=&-2\tilde{k}\alpha^{-2\tilde{k}-1}c(\mu)^{-1/\alpha}P_{\nu-e_0}+\alpha^{-2\tilde{k}-2}c(\mu)^{-1/\alpha}\log c(\mu)P_{\nu-e_0}+\alpha^{-2\tilde{k}}c(\mu)^{-1/\alpha}\partial_{x_2}P_{\nu-e_0}=\nonumber\\
			&=&\alpha^{-2\tilde{k}-2}c(\mu)^{-1/\alpha}\left(-2\tilde{k}\alpha P_{\nu-e_0}+\log c(\mu)P_{\nu-e_0}+\alpha^2\partial_{x_2}P_{\nu-e_0}\right).\nonumber
	\end{eqnarray}}
	Then, the claim holds.
	
	\textbf{3. }$\nu_{i}\geqslant 1$ for some $0<i<N+1$. In this case, $\partial_{\nu}s(\alpha;\mu,z)=\partial_{\mu_i}\partial_{\nu-e_i}s(\alpha;\mu,z)$, where $e_i\in\mathbb{Z}_{\geqslant 0}^{N+2}$ has $1$ at the $(i+1)$th position and zeros elsewhere. As in the previous case, since $|\nu-e_i|=\tilde{k}$, we have
{\small\begin{eqnarray}
			\partial_{\mu_i}\partial_{\nu-e_i}s(\alpha;\mu,z)&=&\partial_{\mu_i}\left(\alpha^{-2\tilde{k}}c(\mu)^{-1/\alpha}P_{\nu-e_i}\left(\log c(\mu),\alpha,z,c^{-1}(\mu),\partial^{\leqslant \tilde{k}}_{\mu}c(\mu)\right)\right)=\nonumber\\
			&=&\alpha^{-2\tilde{k}}\left(-\frac{1}{\alpha}c(\mu)^{-1/\alpha-1}\frac{\partial c}{\partial \mu_i}(\mu)P_{\nu-e_i}+c(\mu)^{-1/\alpha}\partial_{\mu_i}P_{\nu-e_i}\left(\log c(\mu),\alpha,z,c^{-1}(\mu),\partial^{\leqslant \tilde{k}}_{\mu}c(\mu)\right)\right)=\nonumber\\
			&=&\alpha^{-2\tilde{k}-2}c(\mu)^{-1/\alpha}\left(-\alpha c(\mu)^{-1}\frac{\partial c}{\partial \mu_i}(\mu)P_{\nu-e_i}+\alpha^2\partial_{\mu_i}P_{\nu-e_i}\left(\log c(\mu),\alpha,z,c^{-1}(\mu),\partial^{\leqslant \tilde{k}}_{\mu}c(\mu)\right)\right).\nonumber
	\end{eqnarray}}
	By Lemma \ref{LemaPoly}, the claims hold true.
	
	Now, since for any $\nu\in\mathbb{Z}_{\geqslant 0}^{N+2}$ with $|\nu|=k\leqslant K$, we have $$\partial_{\nu}s(\alpha;\mu,z)=\alpha^{-2k}c(\mu)^{-1/\alpha}P_\nu\left(\log c(\mu),\alpha,z,c^{-1}(\mu),\partial^{\leqslant k}_{\mu}c(\mu)\right).$$
	Since $\lim\limits_{\alpha\to 0^{+}}\alpha^{-L}c(\mu)^{-1/\alpha}=0$ for any $L>0$ and $P_\nu$ is a real polynomial and thus is limited for $(\alpha;\mu,z)$ in a neighborhood of $(0;\mu_0,0)$ then $\lim\limits_{\alpha\to 0^{+}}\alpha^{-L}|\partial_{\nu}s(\alpha;\mu,z)|=0$ uniformly and thus all the derivatives of order less or equal to $K$ are bounded, which implies that $s(\alpha;\mu,z)\in\mathcal{F}_{L}^{K}(\mu_0,0)$.
\end{prova}

\begin{lem}\label{LemaSFlatC}
	Let $h\in\mathscr{C}^1(0,\delta)$, for $\delta>0$, such that $h(\alpha)>0$ for $\alpha>0$,
	$$\lim\limits_{\alpha\to 0^{+}}h(\alpha)=\lim\limits_{\alpha\to 0^{+}}\alpha h'(\alpha)=\lim\limits_{\alpha\to 0^{+}}\frac{\alpha\log\alpha}{h(\alpha)}=0.$$
	Then, for any $L>0$, $s(\alpha;z)=\left(1+h(\alpha)\right)^{-1/\alpha}(1+z)\in\mathcal{F}_L^{1}(0)$.
\end{lem}
\begin{prova}
	We have that
$$\log\left(\dfrac{s(\alpha;z)}{\alpha^L}\right)=\log\left(\frac{(1+h(\alpha))^{-1/\alpha}(1+z)}{\alpha^L}\right)=-\frac{1}{\alpha}\log(1+h(\alpha))-L\log\alpha+\log(1+z).$$
The hypothesis $\lim\limits_{\alpha\to 0^{+}}\frac{\alpha\log\alpha}{h(\alpha)}=0$ implies that $\lim\limits_{\alpha\to 0^{+}}\frac{h(\alpha)}{\alpha}=+\infty$. Thus, 
	$$\lim\limits_{\alpha\to 0^{+}}\log (\alpha^{-L}s(\alpha;z))=\lim\limits_{\alpha\to 0^{+}}-\frac{h(\alpha)}{\alpha}\left(1-L\frac{\alpha\log\alpha}{h(\alpha)}\right)+\log(1+z)=-\infty,$$
	which implies that $\lim\limits_{\alpha\to 0^{+}}\alpha^{-L}s(\alpha;z)=0$ and therefore $s(\alpha;z)\in\mathcal{F}_{L}^{0}(\mathbb{R})\subset \mathcal{F}_{L}^{0}(0)$. Now taking the first derivatives of $s(\alpha;z)$ we obtain:
	\begin{eqnarray}
		\partial_z s(\alpha;z)&=&\left(1+h(\alpha)\right)^{-1/\alpha}=\frac{1}{1+z}s(\alpha;z)\nonumber\\
		\partial_\alpha s(\alpha;z)&=&\frac{s(\alpha;z)}{\alpha^{2}}\left(\log(1+h(\alpha))-\alpha\frac{h'(\alpha)}{1+h(\alpha)}\right),\nonumber
	\end{eqnarray}
	since $\lim\limits_{\alpha\to 0^{+}}\alpha^{-L}s(\alpha;z)=0$ and, by hypothesis, $\lim\limits_{\alpha\to 0^{+}}\alpha h'(\alpha)=0$, we have that 
	$$\lim\limits_{\alpha\to 0^{+}}\alpha^{-L}\partial_\alpha s(\alpha;z)=\lim\limits_{\alpha\to 0^{+}}\alpha^{-L}\partial_z s(\alpha;z)=0.$$ 
	Hence, $s(\alpha;z)\in\mathcal{F}_{L}^{1}(\mathbb{R})\subset\mathcal{F}_{L}^{1}(0)$.
\end{prova}

\begin{lem}\label{LemaSol}
	Let $\mu_0$ be a point in an open subset $\Lambda\subset\mathbb{R}^N$. For a given $k\geqslant 1$, consider the equation
	\begin{equation}\label{eqDs}
		s^{-\alpha}-c(\mu)+f(s;\mu)=0,
	\end{equation}
	where $\alpha\in\mathbb{R}\setminus\{0\}$, $f(s;\mu)\in\mathcal{F}^{k}_{\epsilon}(\mu_0)$, for some $0<\epsilon$ and $c(\mu)$ is a $\mathscr C^k$ function at $\mu=\mu_0$ such that $c(\mu_0)\neq 1$ and $c(\mu_0)>0$. Then, there exists an open neighborhood $\tilde{U}$ of $\mu_0$, an open real interval $I$ with $0$ as an endpoint and a $\mathscr C^k$ function $z(\alpha;\mu)$ defined in a neighborhood of $(0,\mu_0)$ such that $z(0;\mu_0)=0$ and $$s=c(\mu)^{-1/\alpha}(1+z(\alpha;\mu)),$$
	is a solution to equation \refc{eqDs} for $(\alpha;\mu)\in I\times\tilde{U}$. More precisely, for $c(\mu_0)>1$, we have $I=(0,\delta)$ and if $c(\mu_0)<1$, then $I=(-\delta,0)$, for a sufficiently small $\delta>0$.
\end{lem}

\begin{prova}
	We have that $s=c(\mu)^{-1/\alpha}$ is a solution of equation \eqref{eqDs} disregarding the $\mathcal{F}_{\epsilon}^{\infty}(\mu_0)$ terms. It is a well defined function, since $c(\mu_0)>0$. We set the following ansatz to estimate the solution of \eqref{eqDs}.
	$$s(\alpha;\mu,z)=c(\mu)^{-1/\alpha}(1+z).$$
	Substituting into \eqref{eqDs} we obtain:
	
	$$c(\mu)\left((1+z)^{-\alpha}-1\right)+h(\alpha;\mu,z)=0,$$
	where $h(\alpha;\mu,z):=f(c(\mu)^{-1/\alpha}(1+z);\mu)$. Note that $\left((1+z)^{-\alpha}-1\right)=-\alpha z g(\alpha,z)$ with $g$ being an analytic function at $(0,0)$ with $g(0,0)=1$. The above equation is equivalent, for $\alpha\neq 0$, to
	$$-c(\mu) z g(\alpha,z)+\tfrac{1}{\alpha}h(\alpha;\mu,z)=0,$$
	Define the functions
	\begin{eqnarray}
		H(\alpha;\mu,z)&=&-c(\mu) z g(\alpha,z)+\tfrac{1}{\alpha}h(\alpha;\mu,z),\nonumber\\
		G(\alpha;\mu,z)&=&\tfrac{1}{\alpha}h(\alpha;\mu,z).\nonumber
	\end{eqnarray}
	Now, since $c(\mu_0)\neq 1$, there is a neighborhood of $\mu_0$ such that $c(\mu)-1$ maintains its sign. Notice that $s(\alpha;\mu,z)$ has different behaviors depending on the sign of $c(\mu)-1$. However, the behaviors are similar regarding the flatness of the function as $\alpha$ goes to zero. More precisely, for $c(\mu)>1$ (and $c(\mu)<1$), $s(\alpha;\mu,z)\to 0$ exponentially flat as $\alpha\to 0^{+}$ (respectively $\alpha\to 0^{-}$).
	
	Suppose that $c(\mu_0)>1$. The case $c(\mu_0)<1$ is similar with the only difference is that the interval $I$ in the statement of the theorem will lie in a different semi-axis. By Lemma \ref{LemaSFlat}, $s(\alpha;\mu,z)\in\mathcal{F}_{L}^{\infty}(\mu_0,0)$ for any $L>0$. Now, since $f(s;\mu)\in\mathcal{F}^\infty_{\epsilon}(\mu_0)$, using item (h) in Lemma \ref{LemaPropFlat}, with $W=(\mu_0)$, $W'=(\mu_0,0)$, $\eta=(\mu,z)$ $\phi(\mu,z)=\mu$, we have that $f(s(\alpha;\mu,z);\mu)\in\mathcal{F}_{\epsilon L}^{\infty}(\mu_0,0)$. Now,  item (g) in Lemma \ref{LemaPropFlat} implies that $G(\alpha;\mu,z)\in\mathcal{F}_{\epsilon L-1}^{\infty}(\mu_0,0)$.
	
	Given a non negative integer $k\geqslant 1$, choosing $L>(k+1)/\epsilon$, by Lemma~\ref{LemaExt}, there exists an extension $\tilde{G}(\alpha;\mu,z)$ of $G(\alpha;\mu,z)$, which is $\mathscr C^k$ on an open neighborhood of $(0;\mu_0,0)$ such that $\tilde{G}(0;\mu_0,0)=0$ and $\partial_z\tilde{G}(0;\mu_0,0)=0$. Consequentially, $H(\alpha;\mu,z)$ admits an extension $\tilde{H}(\alpha;\mu,z)$ which is $\mathscr C^k$ on the same open neighborhood of $(0;\mu_0,0)$ such that
	$$\tilde{H}(0;\mu_0,0)=0,\quad \partial_z\tilde{H}(0;\mu_0,0)=c(\mu_0)\neq 0.$$
	By the Implicit Function Theorem, there exists a unique $\mathscr C^k$ function $z(\alpha;\mu)$ defined in an open set $(-\delta,\delta)\times\tilde{U}$ containing $(0;\mu_0)$ for which
	$$\tilde{H}(\alpha;\mu,z(\alpha,\mu))\equiv 0,\quad \text{and } z(0;\mu_0)=0.$$
	Since $H$ is the restriction of $\tilde{H}$, $s=s(\alpha;\mu,z(\alpha,\mu))$ is a solution for \eqref{eqDs} for $(\alpha;\mu)$ in the open set $I\times\tilde{U}$ for $I=(0,\delta)$ where $\delta>0$ is sufficiently small.
	
	For the case $c(\mu_0)<1$, it is enough to write $s=\hat{c}(\mu)^{-1/\hat{\alpha}}(1+z(\hat{\alpha},\mu))$ with $\hat{c}=1/c$ and $\hat{\alpha}=-\alpha$ and the proof will follow with the same arguments.
\end{prova}

Our goal is to study the zeros of the function
\begin{equation}
	\mathscr{D}(s;\mu)=A(\mu)-B(\mu)s^{-\alpha(\mu)}+\mathcal{F}^{k}_{\epsilon}(\mu_0),
\end{equation}
for $\epsilon>0$, $k\geqslant 1$, with $A(\mu),B(\mu)>0$, $\alpha(\mu)$ being $\mathscr C^k$ functions on $\Lambda$. For any parameter $\mu_0\in\Lambda$, we define $Z(\mathscr{D},\mu_0):=\inf\limits_{\delta,\rho>0} N(\delta,\rho)$, where
$$N(\delta,\rho):=\sup\limits_{\mu\in B_\rho(\mu_0)}\#\{\text{isolated zeros of }\mathscr{D}(\cdot;\mu) \text{ in }(0,\delta)\}.$$
These definitions are motivated by the notion of cyclicity. In fact, the isolated zeros of $\mathscr{D}(s;\mu)$ in \eqref{eqDsAB} correspond to limit cycles of a vector field near a persistent polycycle and, in that case, $Z(\mathscr{D},\mu_0)$ is precisely the cyclicity of the polycycle at parameter value $\mu=\mu_0$. 

\begin{obs}\label{obsrneq1}
For $\alpha(\mu_0)\neq 0$, we must have $Z(\mathscr{D},\mu_0)=0$. In fact, for $\alpha(\mu_0)<0$, we have:
$$\lim\limits_{(s,\mu)\to(0,\mu_0)}\mathscr{D}(s;\mu)=A(\mu_0)>0.$$
And for $\alpha(\mu_0)>0$, 
$$\lim\limits_{(s,\mu)\to(0,\mu_0)}s^{\alpha}\mathscr{D}(s;\mu)=-B(\mu_0)<0.$$
In both cases, there exist an open neighborhood $U$ of $\mu_0$ and $\varepsilon>0$ small enough such that $\mathscr{D}(s;\mu)\neq 0$ for $(s;\mu)\in(0,\varepsilon)\times U$. 
This implies that no limit cycle bifurcate from a persistent polycycle with graphic number different from one.
\end{obs}

Thus, for zeros of $\mathscr{D}(\cdot;\mu)$ to exist near $s=0$ for a small neighborhood of $\mu_0$, we must have $\alpha(\mu_0)=0$.
For the proof of the next theorem, we need to define the \emph{\'Ecalle--Roussarie compensator}.
\begin{defi}
	The function defined for $s>0$ and $\alpha\in\mathbb{R}$ by means of
	$$\omega(s;\alpha)=\left\{\begin{array}{c}
		\frac{s^{-\alpha}-1}{\alpha}\quad\text{if }\alpha\neq 0,\\
		-\log s \quad\text{if }\alpha=0,
	\end{array}\right.$$
	is called \emph{\'Ecalle--Roussarie compensator}.
\end{defi}

\begin{obs}\label{obs-omega}
It is known from  \cite[Lemma A.3]{MarVilDulacGeneral} that $\omega(s;\alpha),\frac{1}{\omega(s;\alpha)}\in\mathcal{F}^\infty_{-\delta}(0)$ for every $\delta>0$.
\end{obs}

\begin{theo}\label{TeoCic1}
	For a fixed value $\mu_0\in\Lambda$, consider the function	
	\begin{equation}\label{eqDsAB}
		\mathscr{D}(s;\mu)=A(\mu)-B(\mu)s^{-\alpha(\mu)}+\mathcal{F}^{k}_{\epsilon}(\mu_0),
	\end{equation}
	where $\epsilon>0$ and $A,B,\alpha\in \mathscr C^k(\Lambda)$ with $k\in\mathbb{N}\cup\{\infty\}$. Assume that
	\begin{itemize}
		\item[$(i)$] $\alpha(\mu)$ changes sign at $\mu_0$, in particular $\alpha(\mu_0)=0$;
		\item[$(ii)$] $A(\mu_0)$ and $B(\mu_0)$ are positive and different.
	\end{itemize}
	Then $Z(\mathscr{D},\mu_0)=1$. Moreover there exist a neighborhood $V$ of $\mu_0$ and $\delta_0>0$ such that $W:=\{\mu\in V:\alpha(\mu)(A-B)(\mu)>0\}\neq\emptyset$ and for every $\mu\in W$ there is an unique $s(\mu)\in(0,\delta_0)$ for which $\mathscr{D}(s(\mu);\mu)=0$. Furthermore, there exists $\hat{z}\in \mathscr C^k(V)$, with $\hat{z}(\mu_0)=0$ such that $s(\mu)=\left(\frac{A}{B}(\mu)\right)^{-1/\alpha(\mu)}\left(1+\hat{z}(\mu)\right)$ for all $\mu\in W$.
\end{theo}

\begin{prova}
	We first show that $Z(\mathscr{D},\mu_0)\leqslant 1$, we express the function $\mathscr{D}(s;\mu)$ as follows
	$$\mathscr{D}(s;\mu)=(A-B)(\mu)-\alpha(\mu)B(\mu)\omega(s;\alpha(\mu))+\mathcal{F}_\epsilon^k(\mu_0),$$
	and apply the derivation-division algorithm. Let $$D_1(s;\mu):=\frac{\mathscr{D}(s;\mu)}{\omega(s;\alpha(\mu))}=-\alpha(\mu)B(\mu)+\frac{(A-B)(\mu)}{\omega(s;\alpha(\mu))}+\mathcal{F}_{\epsilon-\delta}^{k}(\mu_0).$$
	Here, we used the fact that $1/\omega(s;\alpha(\mu))\in\mathcal{F}_{-\delta}^{\infty}(\mu_0)$ for any $\delta>0$ (Remark~\ref{obs-omega}).
	 Then,
	$$\partial_sD_1(s;\mu)=\frac{(A-B)(\mu)s^{-\alpha-1}}{\omega^2(s;\alpha(\mu))}+\mathcal{F}_{\epsilon-\delta-1}^{k}(\mu_0).$$
Since
$\omega(s;\alpha),s^\alpha\in\mathcal{F}^k_{-\delta}(\alpha=0)$, taking $\delta\in(0,\epsilon/4)$, we obtain
$$\lim\limits_{(s,\mu)\to(0,\mu_0)}s^{\alpha+1}\omega^2(s;\alpha(\mu))\partial_sD_1(s;\mu)=\lim\limits_{(s,\mu)\to(0,\mu_0)}(A-B)(\mu)+\mathcal{F}^k_{\epsilon-4\delta}(\mu_0)=(A-B)(\mu_0)\neq 0,$$
	which implies, using Rolle's Theorem, that there are small enough neighborhood $U$ of $\mu_0$ and $\varepsilon>0$ such that $D_1(\cdot;\mu)$ and so $\mathscr{D}(\cdot;\mu)$ have at most one zero for $0<s<\varepsilon$, i.e. $Z(\mathscr{D},\mu_0)\leqslant 1$.
	
	Now, to show that $Z(\mathscr{D},\mu_0)=1$, we use the fact that these zeros are solutions of the following equation: 
	\begin{equation}\label{eqDemoTeo1}
		s^{-\alpha(\mu)}-\dfrac{A(\mu)}{B(\mu)}+\mathcal{F}_{\epsilon}^{k}(\mu_0)=0.
	\end{equation}
	Since $(A-B)(\mu_0)\neq 0$, by Lemma \ref{LemaSol}, there is a solution of the form $s=\left(\frac{A}{B}(\mu)\right)^{-1/\alpha}\left(1+z(\alpha;\mu)\right)$ for $(\alpha;\mu)$ in an open set $I\times\tilde{U}\subset \mathbb{R}\times\mathbb{R}^N$ for equation \eqref{eqDemoTeo1} if we consider $\alpha$ as a free variable. However, since $\alpha(\mu)$ changes signs at $\mu_0$, there exists a neighborhood $V$ of $\mu_0$ such that $W:=\{\mu\in V:\alpha(\mu)(A-B)(\mu)>0\}\subset\tilde{U}\cap\alpha^{-1}(I)$ is not empty. Hence, for every $\mu\in W$, $s(\mu)=\left(\frac{A}{B}(\mu)\right)^{-1/\alpha(\mu)}\left(1+z(\alpha(\mu);\mu)\right)$ is such that $\mathscr{D}(s(\mu);\mu)=0$. Therefore, $Z(\mathscr{D},\mu_0)=1$.
\end{prova}

\section{Proof of Theorem \ref{TeoPeriodmPoly}}\label{SecTeoA}

In this section, we use the results above to describe the behavior of the period function of the limit cycle unfolded from a persistent polycycle. More precisely, we will consider the case in which the cyclicity of this polycycle is exactly 1 at a parameter value $\mu=\mu_0$.

For a family $\{X_\mu\}_{\mu\in\Lambda}$ of polynomial vector fields, we may find an unbounded polycycle $\Gamma$. Thus, to investigate the behavior of the trajectories a given vector field $X_\mu$ near infinity we can consider its Poincaré compactification $p(X_\mu)$ (see \cite[\S 5]{Amarelin}, for details), which is an analytically equivalent vector field defined on $\mathbb{S}^2$. The points at infinity of $\mathbb{R}^2$ are in bijective correspondence with the points at the equator of $\mathbb{R}^2$, denoted by~$\ell_\infty$. Furthermore the trajectories of $p(X_\mu)$ in $\mathbb{S}^2$ are symmetric with respect to the origin and so it suffices to draw only its flow in the closed northern hemisphere, the so-called Poincaré disc. 

\begin{figure}[t]
	\begin{center}		
		\begin{overpic}[width=0.35\textwidth]{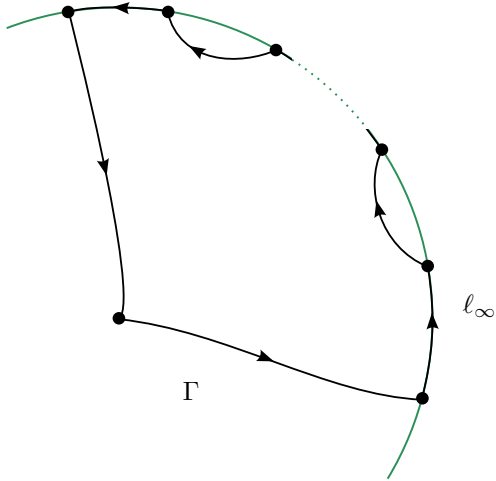}
			\put(37,17) {$\Gamma$}
			\put(95,35) {$\ell_\infty$}
		\end{overpic}
	\end{center}
	\caption{A representation of a polycycle that suits our setting.}\label{Figmcycle}
\end{figure}

For us to have access to the entirety of the previous results, we need to work with a family $\{X_\mu\}_{\mu\in\Lambda}$ of polynomial vector fields such that for $\mu=\mu_0$, $\Gamma$ is a persistent polycycle for $\{X_\mu\}_{\mu\in\Lambda}$. If $\Gamma\cap\ell_\infty\neq\emptyset$, we assume additionally that the infinite line $\ell_\infty$ is invariant for $p(X_\mu)$ for all $\mu\in\Lambda$ (see Figure \ref{Figmcycle}).

The return map $\mathscr{R}(\cdot;\mu)$ associated to a hyperbolic polycycle is broadly studied in the literature. Its general expression is $\mathscr{R}(s;\mu)=A(\mu)s^{r(\mu)}+o(s^{r(\mu)})$ (see \cite{GasMan2002,Soto81,MarVilKolmogorov}, for instance) where $r(\mu)$ is the product of the hyperbolicity ratios of all the saddles in the polycycle, also referred to as \emph{graphic number}. It is known that if $r(\mu_0)\neq 1$ then the cyclicity of the polycycle at $\mu=\mu_0$ is zero and if $r(\mu_0)=1$ and $A(\mu_0)\neq 1$ then the cyclicity is at most 1. Since we are interested in the case which the cyclicity of the persistent polycycle is exactly 1, we work with the additional assumption that $A(\mu)$ changes sign at $\mu=\mu_0$. We now prove Theorem \ref{TeoPeriodmPoly}.

\begin{prooftext}{Proof of Theorem \ref{TeoPeriodmPoly}.}
	For our considered persistent hyperbolic polycycle $\Gamma$, let $p_1,\dots,p_m$ be the $m$ hyperbolic saddles that we conveniently label according to the direction of the flow. We denote by $\lambda_i(\mu)$ the hyperbolicity ratio of $p_i$. First, we consider the case $m\geqslant 2$ (we postpone the proof for the case $m=1$ as it has a slightly different construction). 
	
	Since the infinite line $\ell_\infty$ is invariant for the flow of $p(X_\mu)$, for each saddle located at infinity we have that exactly one of its separatrices is contained in $\ell_\infty$. This also implies that $\Gamma$ has an even number of singularities at infinity assembled in pairs. For $i=1,\dots,m$, let $\Sigma_i$ be the transversal section to the connection from $p_{i-1}$ to $p_i$ (set $p_0=p_m$), and $D_i=D_i(\cdot;\mu)$ and $T_i=T_i(\cdot;\mu)$ the corresponding Dulac maps and times from $\Sigma_i$ to $\Sigma_{i+1}$ (see Figure \ref{Figmcyclesec} for a schematic). It is clear that we can assume that $\Sigma_1$ is the transversal section given in the statement.

	\begin{figure}[t]
		\begin{center}		
			\begin{overpic}[width=0.35\textwidth]{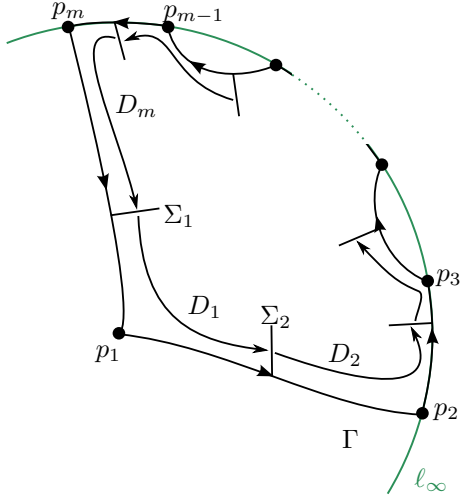}
				\put(70,10) {$\Gamma$}
				\put(19,29) {$p_1$}
				\put(89,17) {$p_2$}
				\put(89.5,45) {$p_3$}
				\put(10,100.5) {$p_m$}
				\put(33,100) {$p_{m-1}$}
				\put(38,38) {$D_1$}
				\put(23,80) {$D_m$}
				\put(67,28) {$D_2$}
				\put(33,57) {$\Sigma_1$}
				\put(53,36) {$\Sigma_2$}
				\put(85,2) {\textcolor{teal}{$\ell_\infty$}}
			\end{overpic}
		\end{center}
		\caption{The representation of the considered maps and saddles on the polycycle $\Gamma$.}\label{Figmcyclesec}
	\end{figure}
	
	By following a point in position $s$ at $\Sigma_1$ until its return to the transversal section, we determine the return map given by
	\begin{equation}\label{eqRispmcycle}
		\mathscr{R}(s;\mu)=D_m\circ D_{m-1}\circ\cdots\circ D_1(s;\mu).
	\end{equation}
	
	Now, in order to study each Dulac map and time, we compactify $X_\mu$ by a means of projective changes of variables instead of the Poincaré compactification $p(X_\mu)$. Note that, for each saddle $p_i$ at infinity, there exists a straight line $l_i$ not intersecting the flow from $\Sigma_i$ to $\Sigma_{i+1}$ for $\mu$ near $\mu_0$ and such that $p_i\not\in l_i\cap l_\infty$. For a fixed $i$, by means of rotation and translation, we can assume that $l_i=\{y=0\}$. Then, we perform the projective change of variables $\{u=\tfrac{x}{y},v=\tfrac{1}{y}\}$ which transforms $X_\mu(x,y)$ into 
	$$\hat{X}^i_\mu(u,v)=\dfrac{1}{v^{d-1}}\left(\hat{P}(u,v;\mu)\partial_u+v\hat{Q}(u,v;\mu)\partial_v\right),$$
	where $d\geqslant 2$ is the degree of $X_\mu$, and $\hat{P},\hat{Q}$ are polynomials with the coefficients depending smoothly on $\mu$.
	It is important to remark that $X_\mu$ and $\hat{X}^i_\mu$ are conjugated, which is essential in the study of the Dulac time. The above transformation brings the infinite line to $\{v=0\}$. This process can be repeated for each saddle at infinity. Now we are ready to work with all saddles in $\Gamma$ (the finite ones, and the infinite ones after the above procedure). For each of them we can apply Lemma \ref{LemaPosarrectes} to conclude that for each $i=1,\dots,m$ there exists a neighborhood $V_i$ of $(p_i,\mu)$ and a $\mathscr C^\infty$ diffeomorphism $\Phi_i:V_i\to\Phi_i(V_i)$ with $\Phi_i(x,y,\mu)=(\phi_\mu^i(x,y),\mu)$ such that $(\phi_\mu^i)_{*}(\hat{X}^i_\mu\vert_{V_i})=\tilde{X}_\mu^i$, where
	\begin{equation}\label{eqXpols}
		\tilde{X}_\mu^i(u,v)=\dfrac{1}{v^{\kappa_i}}\left(uP_i(u,v;\mu)\partial_u+vQ_i(u,v;\mu)\partial_v\right) \mbox{with } \kappa_i=\left\{\begin{array}{l}
			d-1, \mbox{if } p_i\in l_\infty,\\
			0, \mbox{if } p_i\not\in l_\infty,
		\end{array}\right.
	\end{equation}
	where $P_i,Q_i\in \mathscr C^\infty(\Phi(V_i))$. Now, we can apply Theorem \ref{TeoTemps} to obtain that each Dulac map is given by
	$$D_i(s;\mu)=s^{\lambda_i}\left(\Delta_{i0}(\mu)+\mathcal{F}_{\epsilon}^{\infty}(\mu_0)\right)=\Delta_{i0}(\mu)s^{\lambda_i}\left(1+\mathcal{F}_{\epsilon}^{\infty}(\mu_0)\right),$$
	where $\Delta_{i0}\in\mathscr{C}^\infty(\Lambda)$ and we choose $\epsilon>0$ small enough so $0<2\epsilon<\min\{1,\lambda_1(\mu_0),\dots,\lambda_m(\mu_0)\}$. This also implies that $D_i\in\mathcal{F}^{\infty}_{\lambda_i(\mu_0)-\epsilon}(\mu_0)\subset \mathcal{F}^{\infty}_{\epsilon}(\mu_0)$. Using the properties listed in Lemma \ref{LemaPropFlat}, we can see that:
	\begin{eqnarray}
		D_2\circ D_1 (s;\mu)&=&\Delta_{20}(\mu)\left(\Delta_{10}(\mu)s^{\lambda_1}\left(1+\mathcal{F}_{\epsilon}^{\infty}(\mu_0)\right)\right)^{\lambda_2}\left(1+\mathcal{F}_{\epsilon^2}^{\infty}(\mu_0)\right)=\nonumber\\
		&=&\Delta_{20}(\mu)\Delta_{10}^{\lambda_2}(\mu)s^{\lambda_1\lambda_2}\left(1+\mathcal{F}_{\epsilon}^{\infty}(\mu_0)\right)^{\lambda_2}\left(1+\mathcal{F}_{\epsilon^2}^{\infty}(\mu_0)\right)=\nonumber\\
		&=&\Delta_{20}(\mu)\Delta_{10}^{\lambda_2}(\mu)s^{\lambda_1\lambda_2}\left(1+\mathcal{F}_{\epsilon}^{\infty}(\mu_0)\right)\left(1+\mathcal{F}_{\epsilon^2}^{\infty}(\mu_0)\right)=\nonumber\\
		&=&\Delta_{20}(\mu)\Delta_{10}^{\lambda_2}(\mu)s^{\lambda_1\lambda_2}\left(1+\mathcal{F}_{\epsilon^2}^{\infty}(\mu_0)\right),\nonumber
	\end{eqnarray}
	where in the third equality we use the fact that $(s;\mu)\to(1+s)^{\lambda_i(\mu)}-1$ belongs to $\mathcal{F}^\infty_{1}(\mu_0)$. Notice that the composition of two Dulac maps has a Dulac map like expression. Proceeding analogously, we obtain the composition of $k$ consecutive Dulac maps:
	\begin{equation}\label{eqkDulacs}
		D_k\circ D_{k-1}\circ\cdots\circ D_1(s;\mu)=A_{k}(\mu)s^{\lambda_1\dots\lambda_k}\left(1+\mathcal{F}_{\epsilon^k}^{\infty}(\mu_0)\right),
	\end{equation}
	where 
	$$A_k(\mu)=\prod_{i=1}^{k}\Delta_{i0}^{\Lambda_{ik}}(\mu),\quad \Lambda_{ik}=\prod_{i<j\leqslant k}\lambda_j \text{ and } \Lambda_{kk}=1.$$
	Therefore, the return map is given by
	\begin{equation}\label{eqReturn}
		\mathscr{R}(s;\mu)=A_m(\mu)s^{\lambda_1\dots\lambda_m}\left(1+\mathcal{F}_{\epsilon^m}^{\infty}(\mu_0)\right),
	\end{equation}
	which proves (a).
	
	The return time of a point in position $s$ at $\Sigma_1$ is determined by the following expression:
	\begin{equation}\label{eqTempsPoly}
		\mathcal{T}(s;\mu)=T_1(s;\mu)+\sum_{i=2}^{m}T_{i}(D_{i-1}\circ\dots\circ D_{1}(s;\mu)).
	\end{equation}
	For each saddle $p_i$ in $\Gamma$, by \eqref{eqXpols}, the vector field can be locally expressed as \eqref{eqX1} with $\mathtt{n}=(0,n_2)$. Therefore, by Theorem \ref{TeoTemps}, the associated time map $T_i$ is given by
	$$T_i(s;\mu)=T_0^i(\mu)\log s+T_{00}^{i}(\mu)+\mathcal{F}_{\epsilon}^{\infty}(\mu_0)=T_0^i(\mu)\log s+\mathcal{F}_{0}^{\infty}(\mu_0).$$
	
	By \eqref{eqkDulacs} and Lemma \ref{LemaPropFlat} for $i\geqslant 2$, we have
	\begin{eqnarray}
		T_{i}(D_{i-1}\circ\dots\circ D_{1}(s;\mu))&=&T_{0}^{i}(\mu)\Lambda_{0 (i-1)}\log s+T_{0}^{i}(\mu)\log A_{i-1}(\mu)+\mathcal{F}_{0}^{\infty}(\mu_0)=\nonumber\\
		&=&T_{0}^{i}(\mu)\Lambda_{0 (i-1)}\log s+\mathcal{F}_{0}^{\infty}(\mu_0).\nonumber
	\end{eqnarray}
	Thus, substituting into \eqref{eqTempsPoly}, we obtain
	\begin{equation}\label{eqTempsPoly1}
		\mathcal{T}(s;\mu)=\bar{T}_0(\mu)\log s +\mathcal{F}^{\infty}_{0}(\mu_0),
	\end{equation}
	for $\bar{T}_{0}(\mu)=T_{0}^1(\mu)+\sum_{i=2}^{m}T_{0}^{i}(\mu)\Lambda_{0 (i-1)}$. Note that, by Theorem \ref{TeoTemps}, $T_{0}^{i}\leqslant 0$ and the equality holds if and only if $p_i$ is an infinite saddle. If there is at least one finite saddle, then $\bar{T}_{0}(\mu)$ is a strictly negative function. This proves (b).
	
	The limit cycles bifurcating from $\Gamma$ correspond to small positive zeros of the proportional displacement map $\mathscr{D}(s;\mu)=s^{-\lambda_1\dots\lambda_n}(\mathscr{R}(s;\mu)-s)$, which by \eqref{eqReturn} writes as
	\begin{equation}\label{eqDispmcycle}
		\mathscr{D}(s;\mu)=A_m(\mu)-s^{1-\lambda_1\dots\lambda_m}+\mathcal{F}_{\epsilon^m}^{\infty}(\mu_0).
	\end{equation}
	Recall that, by assumption, $r(\mu)-1=(\lambda_1\dots\lambda_m)(\mu)-1$ changes sign at $\mu_0$ and $A_m(\mu_0)\neq 1$. Hence, by Theorem \ref{TeoCic1},  $\Gamma$ has cyclicity 1 at $\mu_0$. Furthermore, there exist a neighbourhood $V$ of $\mu_0$ and a smooth function $z\in\mathscr C^\infty(V)$ with $z(\mu_0)=0$ such that the location of the bifurcating limit cycle is given by $s=s(\mu)$ with
	\begin{equation}\label{eqSPoly}
		s(\mu)=A_m(\mu)^{-1/(r(\mu)-1)}\big(1+z(\mu)\big),
	\end{equation}
	for $\mu\in W:=\{\mu\in V:(r(\mu)-1)(A_m(\mu)-1)>0\}$.
	Then, choosing $K_0\in(0,|\log(A_m(\mu_0))|)$ we have
	\begin{equation}\label{eqExpfast}
		\lim\limits_{\substack{\mu\to\mu_0\\ \mu\in W}}\dfrac{s(\mu)}{\exp\left(\frac{-K_0}{|1-r(\mu)|}\right)}=\lim\limits_{\substack{\mu\to\mu_0\\ \mu\in W}}\exp(|\log(A_m(\mu))|-K_0)^{\frac{-1}{|1-r(\mu)|}}(1+z(\mu))=0,
	\end{equation}
	which concludes the proof of (c). With regard to the case $m\geqslant 2$ it only remains to be proved the assertion (d) about the period of the bifurcating limit cycle. To this end we assume that at least one of the saddles in $\Gamma$ is finite. The period $\mathscr{T}(\mu)$ of the limit cycle $\gamma_\mu$ bifurcating from $\Gamma$ is obtained inserting the location of the limit cycle given by \eqref{eqSPoly} into the return time \eqref{eqTempsPoly1}. Hence,
$$\lim\limits_{\substack{\mu\to\mu_0\\ \mu\in W}}|1-r(\mu)|\mathscr{T}(\mu)=\lim\limits_{\substack{\mu\to\mu_0\\ \mu\in W}}|1-r(\mu)|\mathcal{T}(s(\mu);\mu)=\lim\limits_{\substack{\mu\to\mu_0\\ \mu\in W}}-\bar{T}_0(\mu)|\log A_m(\mu)|=-\bar{T}_0(\mu_0)|\log A_m(\mu_0)|> 0.$$
	
	\textbf{Case $m=1$:} To complete the proof, we consider the case where $m=1$, i.e. the case where $\Gamma$ is a (finite) persistent saddle loop. 
	
	\begin{figure}[t]
		\begin{center}		
			\begin{overpic}[width=0.35\textwidth]{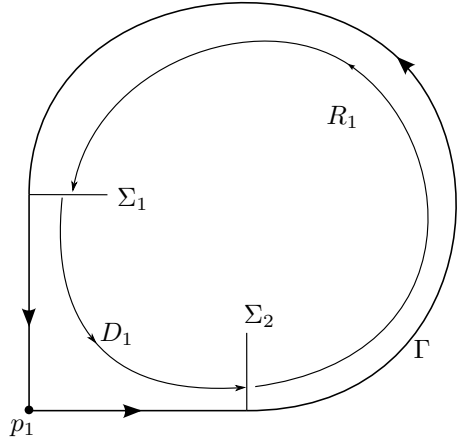}
				\put(-2.5,-3) {$p_1$}
				\put(22,49) {$\Sigma_1$}
				\put(18,18) {$D_1$}
				\put(51,22) {$\Sigma_2$}
				\put(70,68) {$R_1$}
				\put(90,14) {$\Gamma$}
			\end{overpic}
		\end{center}
		\caption{The transversal sections and the maps defined around a homoclinic loop $\Gamma$ of a hyperbolic saddle.}\label{FigSaddleLoop}
	\end{figure}
	
	We place two transversal sections $\Sigma_1$ and $\Sigma_2$ along the homoclinic connection labeled according to the direction of the flow. Now, consider $D_1$ and $T_{1}$ the Dulac map and time from $\Sigma_1$ to $\Sigma_2$ and $R_1$ and $T_{2}$ the corresponding maps associated to the passage from $\Sigma_1$ to $\Sigma_2$ using the flow of $-X_{\mu}$, as depicted in Figure~\ref{FigSaddleLoop}. More precisely, $\sigma_2(R_1(s;\mu))=\varphi(-T_{2}(s;\mu),\sigma_1(s);\mu)$, where $\sigma_i$ is a parametrization of the transversal section $\Sigma_i$, $\varphi(t,q;\mu)$ is the orbit of $X_\mu$ passing through $q\in\mathbb{R}^2$ at $t=0$.
	
	The return map is given by
	$$\mathscr{R}(s;\mu)=R_1\circ D_1(s;\mu).$$
	Note that the function $R_1(s;\mu)$ is a diffeomorphism. Since, the saddle loop is unbroken for $\mu\in\Lambda$, we have that $R_1(0;\mu)=0$. Thus, we obtain that $R_1(s;\mu)=a_1(\mu)s+\mathcal{F}^{\infty}_{2}(\mu_0)$, with $a_1(\mu)$ being a strictly positive function. Exactly as before, by Lemma \ref{LemaPosarrectes} and Theorem \ref{TeoTemps}, we have that $D_1(s;\mu)=s^{\lambda_1}(\Delta_{10}(\mu)+\mathcal{F}_{\epsilon}^{\infty}(\mu_0))$ for $\epsilon\in (0,\min\{\lambda_1(\mu_0),1\})$. Therefore, the return map is given by
	\begin{equation}\label{eqRPeix}
		\mathscr{R}(s;\mu)=\bar{A}(\mu)s^{\lambda_1}\left(1+\mathcal{F}_{\epsilon}^{\infty}(\mu_0)\right),
	\end{equation}
	with $\bar{A}(\mu)=a_1(\mu)\Delta_{10}(\mu)$, which proves (a). The return time $\mathcal{T}(s;\mu)$ of a point in position $s$ at $\Sigma_1$ is:
	$$\mathcal{T}(s;\mu)=T_1(s;\mu)+T_2(D_1(s;\mu);\mu).$$
	Now $T_2(s;\mu)$ is the regular time from $\Sigma_2$ to $\Sigma_1$, which is a $\mathscr{C}^{\infty}$ function in a neighborhood of $(s;\mu)=(0;\mu_0)$. Thus $T_2\in\mathcal{F}_{0}^{\infty}(\mu_0)$. As for the Dulac time $T_1(s;\mu)$, we can argue as in the previous case, using  Lemma \ref{LemaPosarrectes} and Theorem \ref{TeoTemps} to conclude that
	$$T_1(s;\mu)=T_0^1(\mu)\log s+T_{00}^{1}(\mu)+\mathcal{F}_{\epsilon}^{\infty}(\mu_0)=T_0^1(\mu)\log s+\mathcal{F}_{0}^{\infty}(\mu_0),$$
	with $T_0^1(\mu)<0$. Therefore, 
	\begin{equation}\label{eqTempsPeix}
		\mathcal{T}(s;\mu)=T_0^1(\mu)\log s+\mathcal{F}_{0}^{\infty}(\mu_0).
	\end{equation}
	This proves (b). By \eqref{eqRPeix}, the proportional displacement map is given by:
	\begin{equation}
		\mathscr{D}(s;\mu)=s^{-\lambda_1}\left(\mathscr{R}(s;\mu)-s\right)=\bar{A}(\mu)-s^{1-\lambda_1}+\mathcal{F}_{\epsilon}^{\infty}(\mu_0).
	\end{equation}
	By Theorem \ref{TeoCic1}, the sufficient conditions for $\Gamma$ to have cyclicity 1 at $\mu_0$ are $r(\mu)-1=\lambda_1(\mu)-1$ changes sign at $\mu_0$ and $\hat{A}(\mu_0)\neq 1$ and under these conditions, the location of the lone limit cycle is given by 
	$$s(\mu)=\bar{A}(\mu)^{-1/(r(\mu)-1)}(1+z(\mu)),$$
	for $\mu\in W:=\{\mu\in V:(r(\mu)-1)(\bar{A}(\mu)-1)>0\}$, $V$ being a small enough neighborhood of $\mu_0$. Here we can perform the same computations as in \eqref{eqExpfast} to prove (c).
	
	The period $\mathscr{T}(\mu)$ of a limit cycle $\gamma_\mu$ unfolded from $\Gamma$ at the parameter value $\mu$ is given by $\mathscr{T}(\mu)=\mathcal{T}(s(\mu);\mu)$. Thus, substituting $s(\mu)$ into \eqref{eqTempsPeix}, we conclude that the asymptotic behavior of the period of the limit cycle $\gamma_\mu$ is
	$$\lim\limits_{\substack{\mu\to\mu_0\\ \mu\in W}}(\lambda_1(\mu)-1)\mathscr{T}(\mu)=
	-T_{0}^{1}(\mu_0)|\log(\bar{A}(\mu_0))|,$$
	which proves (d).
\end{prooftext}

\begin{obs}
	If $A_m(\mu_0)>1$, the limit must be taken as $r(\mu)\to 1^{+}$ and for $A_m(\mu_0)<1$,  $r(\mu)\to 1^{-}$. This assures we are always at a well defined position $s(\mu)$ for the limit cycle $\gamma_\mu$ as we approach the polycycle~ $\Gamma$.
\end{obs}

Theorem \ref{TeoPeriodmPoly} provides sufficient conditions for the cyclicity of the persistent polycycle to be exactly 1 and also describes how fast the unfolded limit cycle approaches the polycycle. Moreover, asymptotic behavior of the period of the limit cycle is the same as the function $\frac{1}{|1-r(\mu)|}$ where $r(\mu)$ is the graphic number. This result encompasses a wide range of polycycles, among them the Kolmogorov polycycles studied in \cite{MarVilKolmogorov}.

\section{Proof of Theorem \ref{TeoCenter}}\label{SecTeoB}

Up to this point, we studied families of polynomial systems with a persistent polycycle whose cyclicity is 1 at a parameter $\mu_0$ but the return map near the polycycle is not the identity. In those families, by Theorem~\ref{TeoPeriodmPoly}, the period of the limit cycle unfolded from the polycycle has a determined behavior. However, that may not be the case when the return map near the polycycle is the identity at $\mu_0$. In this case, we say that the polycycle is of the \emph{center-type}.

In \cite[Theorem 3.1, item (c)]{MarVilKolmogorov}, the authors exhibit an example of a center-type polycycle. It is present in the following family
\begin{equation}\label{eqKolmogorovC}
	X_{\mu}=x(1+x+x^2+axy+py^2)\partial_x+y(-1-y+qy^2+axy-y^2)\partial_y,
\end{equation}
where $\mu=(a,p,q)\in\mathbb{R}^3$ with $p<-1$, $q>1$. The displacement function associated to the polycycle of this Kolmogorov system is, after a change of variables in the parameter space, proportional to the following expression
$$\mathscr{D}(s;\mu)=s^{-\alpha(\mu)}-(1+\varepsilon(\mu))+f(s;\mu).$$
Moreover, if $\mu_0=(0,p_0,-p_0)$ then $\mathscr{D}(s;\mu_0)\equiv 0$, $f(s;\mu)\in\mathcal{F}_{\epsilon}^{\infty}(\mu_0)$, $\epsilon>0$ and $\alpha(\mu),\varepsilon(\mu)$ are such that $\alpha(\mu_0)=\varepsilon(\mu_0)=0$ and $\rm{rank}\{\nabla\alpha,\nabla\varepsilon\}\vert_{\mu=\mu_0}=2$.

Note that in this case, we cannot apply Theorem \ref{TeoCic1}. However, the cyclicity of the polycycle at $\mu_0$ is 1.

This example motivates \teoc{TeoCenter}, which will focus in families of polynomial vector fields $\{X_{\mu}\}_{\mu\in\Lambda}$ such that $X_{\mu_0}$ has a persistent hyperbolic polycycle $\Gamma$ with the infinite line $\ell_\infty$ being invariant for all $\mu\in\Lambda$ and with at least one finite saddle, with the additional assumption that for a fixed transversal section we can write the proportional displacement function near the the polycycle as
\begin{equation}
	\mathscr{D}(s;\mu)=s^{-\alpha(\mu)}-(1+\varepsilon(\mu))+f(s;\mu),
\end{equation}
with $f(s;\mu)\in\mathcal{F}_{\epsilon}^{\infty}(\mu_0)$, $\epsilon>0$, with $\alpha(\mu_0)=\varepsilon(\mu_0)=0$ and $\nabla\alpha,\nabla\varepsilon$ are linearly independent at $\mu_0$. These assumptions encompass the cases for which the polycycle $\Gamma$ is of the center-type or the cyclicity of the polycycle is greater than 1. 

\begin{prooftext}{Proof of Theorem \ref{TeoCenter}.}
	By item (a) in Theorem \ref{TeoPeriodmPoly}, the return map associated to $\Gamma$ writes as
	$$\mathscr{R}(s;\mu)=s^{r(\mu)}\left(A(\mu)+f(s;\mu)\right),$$
	with $f(s;\mu)\in\mathcal{F}^{\infty}_{\epsilon}(\mu_0)$ for $\epsilon>0$ small enough. Let $\alpha(\mu)=r(\mu)-1$ and $\varepsilon(\mu)=A(\mu)-1$. Since ${\rm rank}\{\nabla A,\nabla r\}\vert_{\mu=\mu_0}=2$, there exists a local diffeomorphism $\varphi(\mu)=(\alpha(\mu),\varepsilon(\mu),\bar{\mu}(\mu))$ defined on a small neighborhood of $\mu_0$ such that $\varphi(\mu_0)=(0,0,\bar{\mu}_0)$. Thus, we can work locally with the parameter set $(\alpha,\varepsilon,\bar{\mu})$. The parameter $\bar{\mu}$ plays a secondary role in our proof and therefore, we will omit it in the rest of the argument to simplify the notation. In our new parameters, the proportional displacement map $\mathscr{D}(s;\alpha,\varepsilon)=-s^{-\alpha-1}\left(\mathscr{R}(s;\varphi^{-1}(\alpha,\varepsilon))-s\right)$ is given by
	\begin{equation}\label{eqDcenter}
		\mathscr{D}(s;\alpha,\varepsilon)=s^{-\alpha}-(1+\varepsilon)+\bar{f}(s;\alpha,\varepsilon),
	\end{equation}
	where $\bar{f}(s;\alpha,\varepsilon)=f(s;\varphi^{-1}(\alpha,\varepsilon))$. We have $\bar{f}(s;\alpha,\varepsilon)\in\mathcal{F}_{\epsilon}^{\infty}(0,0)$.
	For the given function $\tau$, consider the function $h(\alpha)=\alpha\tau(\alpha)$ and the $\mathscr{C}^1$ curve $m(\alpha)=\varphi^{-1}(\alpha,h(\alpha))$  in the parameter space defined for small $\alpha>0$  (here, we consider $\bar{\mu}=\bar{\mu}_0$ fixed). 
	It is clear that $\lim\limits_{\alpha\to 0^{+}}m(\alpha)=\mu_0$. For parameters on this curve, the proportional displacement map becomes
	\begin{equation}\label{eqDcenter2}
		M(s;\alpha):=\mathscr{D}(s;\alpha,h(\alpha))=s^{-\alpha}-(1+h(\alpha))+\bar{f}(s;\alpha,h(\alpha)).
	\end{equation}
	The zeros of $M(\cdot;\alpha)$ correspond to limit cycles of the corresponding vector field $X_{m(\alpha)}$. To search for these zeros, we establish the following ansatz
	$$s(\alpha;z)=\left(1+h(\alpha)\right)^{-1/\alpha}(1+z).$$
	Substituting into \eqref{eqDcenter2}, we obtain
	\begin{eqnarray}
			H(\alpha,z)&:=&\tfrac{1}{\alpha}M(s(\alpha;z),\alpha)=\tfrac{1}{\alpha}\left(\left(1+h(\alpha)\right)(1+z)^{-\alpha}-(1+h(\alpha))+\bar{f}(s(\alpha;z);\alpha,h(\alpha))\right)=\nonumber\\
			&=&(1+h(\alpha))\omega(1+z;\alpha)+\tfrac{1}{\alpha}\bar{f}(\left(1+h(\alpha)\right)^{-1/\alpha}(1+z);\alpha,h(\alpha)).\nonumber
	\end{eqnarray}
	Since $\omega(1+z;\alpha)=F(\alpha\log(1+z))\log(1+z)$ where $F(x)=\frac{e^{-x}-1}{x}$, the function $(1+h(\alpha))\omega(1+z;\alpha)$ admits a continuous extension $g_1(\alpha,z)$ in a neighborhood of $(\alpha,z)=(0,0)$ such that $g_1(0,0)=0$ and $\partial_zg_1$ is continuous with
	$$\partial_zg_1(\alpha,z)=\tfrac{1+h(\alpha)}{1+z}\left(F'(\alpha\log(1+z))\alpha\log(1+z)+F(\alpha\log(1+z))\right),$$
	and $\partial_zg_1(0,0)=-1$. Moreover, the assumptions on $\tau$ imply that $h$ is a $\mathscr{C}^1$ function for small $\alpha>0$, such that $h(\alpha)>0$ for $\alpha>0$, and $\lim\limits_{\alpha\to 0^{+}}h(\alpha)=\lim\limits_{\alpha\to 0^{+}}\alpha h'(\alpha)=\lim\limits_{\alpha\to 0^{+}}\frac{\alpha\log\alpha}{h(\alpha)}=0$. Thus, by Lemma \ref{LemaSFlatC}, for any $L>0$, $s(\alpha;z)\in\mathcal{F}_L^{1}(0)$. Now, since $\bar{f}(s;\alpha,\varepsilon)\in\mathcal{F}_{\epsilon}^{\infty}(0,0)$, by items (g) and (h) in Lemma \ref{LemaPropFlat}, we have that $g(\alpha;\varepsilon,z)=\alpha^{-1}\bar{f}(s(\alpha;z);\alpha,\varepsilon)\in\mathcal{F}_{\epsilon L-1}^{1}(0,0)$. Taking $L>2/\varepsilon$, by Lemma~\ref{LemaExt}, there is an extension $\tilde{g}(\alpha;\varepsilon,z)$ of $g(\alpha;\varepsilon,z)$ which is $\mathscr C^1$ on an open neighborhood of $(0;0,0)$ such that $\tilde{g}(0;0,0)=\partial_z\tilde{g}(0;0,0)=0$. Thus, we can use this extension to conclude that the function $g(\alpha;h(\alpha),z)$ admits a continuous extension $g_2(\alpha,z)$ with continuous partial derivative $\partial_zg_2$ such that $g_2(0,0)=\partial_zg_2(0,0)=0$.
	
	Now, we have that $H(\alpha,z)$ admits $\tilde{H}(\alpha,z):=g_1(\alpha,z)+g_2(\alpha,z)$ a a continuous extension with continuous partial derivative $\partial_z\tilde{H}$ such that $\tilde{H}(0,0)=0$ and $\partial_z(0,0)=-1$. 
	By the Continuous Implicit Function Theorem \cite[Theorem 2]{Implicit}, there exists a unique function $\hat{z}(\alpha)$, continuous, defined in an open set $(-\delta,\delta)$, for which
	$$\tilde{H}(\alpha;\hat{z}(\alpha))\equiv 0,\;\text{and } \hat{z}(0)=0.$$
	Therefore $s=s(\alpha;\hat{z}(\alpha))$ is a zero for \eqref{eqDcenter2} for $\alpha>0$ small enough. We recall that $s(\alpha;\hat{z}(\alpha))$ corresponds to the location of the intersection of a limit cycle for vector field $X_{m(\alpha)}$ and the transversal section $\Sigma$. Since $\lim\limits_{\alpha\to 0^{+}}s(\alpha;\hat{z}(\alpha))=0$, this limit cycle approaches the polycycle as $\alpha\to 0^+$. We now proceed to compute its period $\mathscr{T}(\alpha)$. By item (b) in Theorem \ref{TeoPeriodmPoly} the return time associated to $\Gamma$ is given by 
	$$\mathcal{T}(s;\mu)=\bar{T}_0(\mu)\log s +\mathcal{F}^{\infty}_{0}(\mu_0),$$
	where $\bar{T}_{0}\in\mathscr{C}^{\infty}(\Lambda)$ with $\bar{T}_{0}(\mu)<0$. Then, we have that $\mathscr{T}(\alpha)=\mathcal{T}(s(\alpha;\hat{z}(\alpha));m(\alpha))$ and
	\begin{eqnarray}
		\lim\limits_{\alpha\to 0{+}}\dfrac{\mathscr{T}(\alpha)}{\tau(\alpha)}&=&\lim\limits_{\alpha\to 0{+}}\frac{1}{\tau(\alpha)}\bar{T}_0(m(\alpha))\log \left(1+h(\alpha)\right)^{-1/\alpha}=\nonumber\\
		&=&\lim\limits_{\alpha\to 0{+}}-\frac{\bar{T}_0(m(\alpha))}{\alpha\tau(\alpha)}\log \left(1+\alpha\tau(\alpha)\right)=-\bar{T}_0(\mu_0)>0.\nonumber
	\end{eqnarray}
This concludes the proof of the result.	
\end{prooftext}

\begin{figure}[t]
	\begin{center}		
		\begin{overpic}[width=0.35\textwidth]{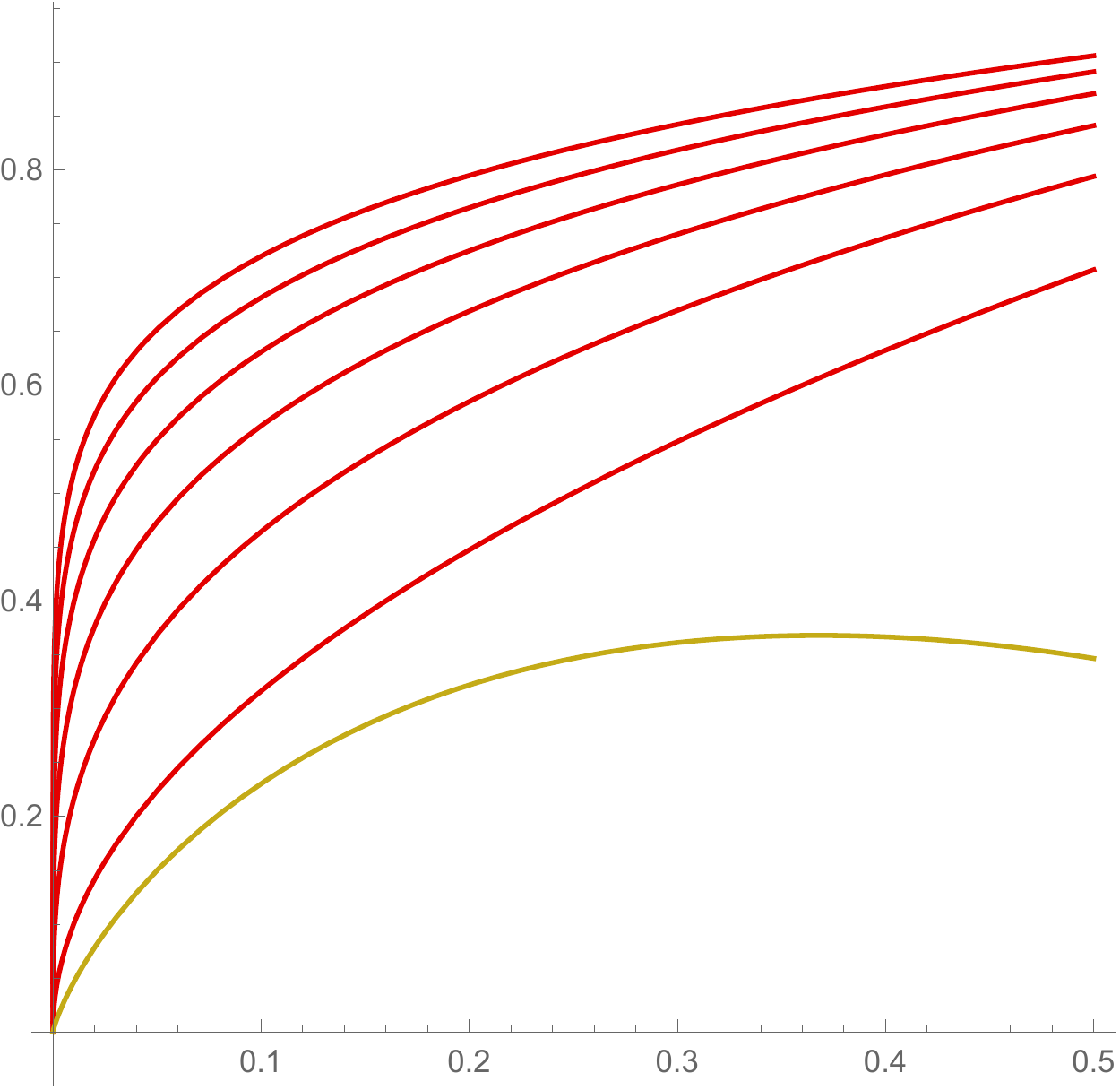}
			\put(101,2) {$\alpha$}
			\put(2,100) {$\varepsilon$}
			\put(101,38) {$-\alpha\log\alpha$}
			\put(95,80) {$\left.\begin{array}{c}\\
					\\
				\end{array}\right\}\alpha^{l+1}$}
		\end{overpic}
	\end{center}
	\caption{The functions $h(\alpha)=\alpha^{l+1}$, for $l\in(-1,0)$ in relation to the threshold $-\alpha\log\alpha$ for which we can apply Lemma~\ref{LemaSFlatC}.}\label{Figh}
\end{figure}
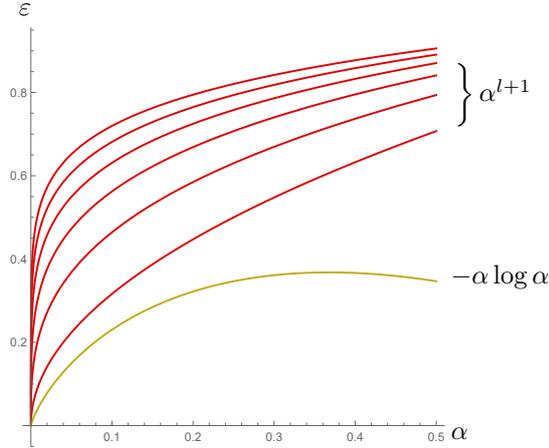

Theorem \ref{TeoCenter} tells us that there is not a single determinate behavior for the period of a limit cycle unfolding from the polycycle for the considered family of vector fields $X_\mu$. In particular, if we consider the functions $\tau(\alpha)=\alpha^l$, for $l\in(-1,0)$, we can always find a path in the parameter space approaching $\alpha=0$ for which the limit cycle has period $\mathscr{T}(\alpha)\sim \alpha^{l}$. This is substantially different from the cases encompassed by Theorem~\ref{TeoPeriodmPoly}, since the behavior of the period of the limit cycle does not depend on the way we approach the parameter value $\mu_0$. In Figure~\ref{Figh} we represent graphically the behavior of the functions $h(\alpha)=\alpha\tau(\alpha)$ compared to the function $-\alpha\log\alpha$ showing some  functions which are suitable to be the period of a limit cycle bifurcating from the polycycle in the case encompassed by Theorem~\ref{TeoCenter}.

\appendix

\section{Dulac map and time}

Since we will deal with persistent hyperbolic polycycles, we will need to work with the Dulac map and Dulac time associated to hyperbolic saddles. In this appendix we defined these concepts in a particular setting where the expression of these maps is known. For more details, we refer the reader to \cite{MarVilDulacLocal,MarVilDulacGeneral} where the specifics are carried out extensively. 

We consider an open set $\Lambda\subset\mathbb{R}^N$ and the family $\{X_{\mu}\}_{\mu\in\Lambda}$ of vector fields given by:
\begin{equation}\label{eqX1}
	X_{\mu}:=\dfrac{1}{x^{n_1}y^{n_2}}\left(xP(x,y;\mu)\partial_x+yQ(x,y;\mu)\partial_y\right).
\end{equation}
Here,
\begin{itemize}
	\item $\mathtt{n}:=(n_1,n_2)\in\mathbb{Z}_{\geqslant 0}^2$;
	\item $P,Q\in \mathscr C^{\infty}(V\times \Lambda),$ for some open set $V\subset\mathbb{R}^2$ containing the origin;
	\item $P(x,0,\mu)>0$ and $Q(0,y,\mu)<0$, for all $(x,0), (0,y)\in V$ and $\mu\in\Lambda$. This means that the origin is a hyperbolic saddle of $x^{n_1}y^{n_2}X_\mu$ with the $y$-axis being the stable manifold and $x$-axis the unstable manifold;
	\item $\lambda(\mu)=-\dfrac{Q(0,0;\mu)}{P(0,0;\mu)}$ is the hyperbolic ratio of the saddle.
\end{itemize}

For $i=1,2$, let $\sigma_i:(-\epsilon,\epsilon)\times\Lambda\to\Sigma_i$ be transverse sections of $X_{\mu}$ to the axis such that

\begin{figure}[t]
	\begin{center}		
		\begin{overpic}[width=0.7\textwidth]{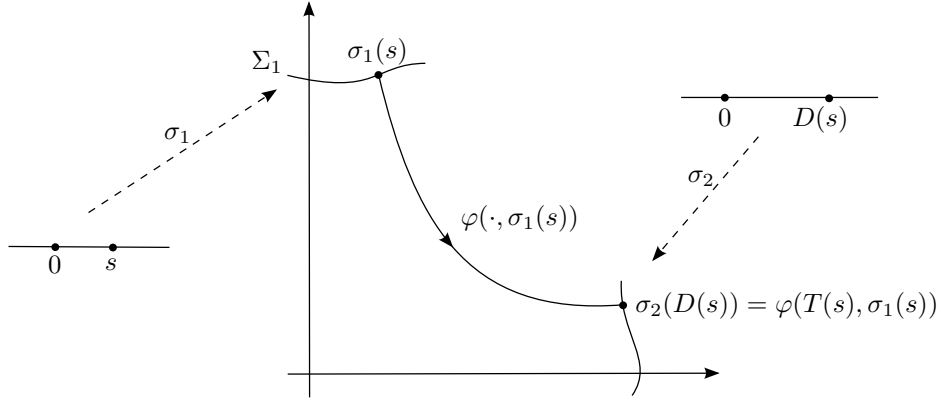}
			\put(4.7,14.5) {$0$}
			\put(11.2,14.8) {$s$}
			\put(18,30) {$\sigma_1$}
			\put(28,38) {$\Sigma_1$}
			\put(52,20) {$\varphi(\cdot,\sigma_1(s))$}
			\put(39,39.3) {$\sigma_1(s)$}
			\put(72,10) {$\sigma_2(D(s))=\varphi(T(s),\sigma_1(s))$}
			\put(78,25) {$\sigma_2$}
			\put(81.5,31.5) {$0$}
			\put(90,31.5) {$D(s)$}
		\end{overpic}
	\end{center}
	\caption{The Dulac map and time.}\label{FigDulac}
\end{figure}

$\sigma_1(0;\mu)\in\{(0,y):y>0\}$ and $\sigma_2(0,\mu)\in\{(x,0):x>0\}$ for all $\mu\in\Lambda$. The Dulac map $D(\cdot;\mu)$ and the Dulac time $T(\cdot;\mu)$ are defined by the following relationship:
$$\varphi(T(s;\mu),\sigma_1(s;\mu);\mu)=\sigma_2(D(s;\mu),\mu),\;\forall s\in (0,\epsilon),$$
where $\varphi(t,p_0;\mu)$ is the orbit of $X_{\mu}$ passing through $p_0\in V$ at $t=0$ (see Figure~\ref{FigDulac}).

The following result gathers Theorems A and B in \cite{MarVilDulacGeneral} to give an asymptotic development of the Dulac map and time.

\begin{theo}\label{TeoTemps}
	Let $D(s;\mu)$ and $T(s;\mu)$ be respectively the Dulac map and time of the hyperbolic saddle~\eqref{eqX1} from $\Sigma_1$ to $\Sigma_2$. Then, there exists $\Delta_{00}\in \mathscr C^\infty(\Lambda)$ such that for any $\lambda_0>0$ and $\epsilon>0$ small enough so $\epsilon<\min\{\lambda_0,1\}$, for which
	$$D(s;\mu)=s^\lambda(\Delta_{00}(\mu)+\mathcal{F}_{\epsilon}^{\infty}(\mu_0)).$$
	And for $\mathtt{n}=(n_1,0)$ or $\mathtt{n}=(0,n_2)$, there exists $T_{00}\in \mathscr C^{\infty}(\Lambda)$ such that for any $\lambda_0>0$ and $\epsilon>0$ small enough so $\epsilon<\min\{\lambda_0,1\}$, we have
	$$T(s;\mu)=T_0(\mu)\log s+T_{00}(\mu)+\mathcal{F}_{\epsilon}^{\infty}(\mu_0),$$
	where $T_0(\mu)=\left\{\begin{array}{cc}
		0, &\text{if } \mathtt{n}\neq (0,0),\\
		\tfrac{-1}{P(0,0;\mu)}, &\text{if }\mathtt{n}=(0,0).
	\end{array}\right.$
\end{theo}

Note that the above result applies to hyperbolic saddles for which the separatrices are contained in the orthogonal axis. However, this is not a restrictive assumption since we can rectify the separatrices of hyperbolic saddles via a family $\mathscr C^\infty$ diffeomorphism \cite[Lemma 4.3]{MarVilDulacGeneral}. For the sake of the reader, we present the result below.

\begin{lem}\label{LemaPosarrectes}
	Consider a $\mathscr C^\infty$ family $\{X_\mu\}_{\mu\in\mathbb{R}^N}$ of planar vector fields defined in an open set of $\mathbb{R}^2$.	For a fixed $\mu_0\in\mathbb{R}^N$ and assume that for all $\mu$ in a neighborhood of $\mu_0$, $X_\mu$ has a hyperbolic saddle $p_\mu$ with stable and unstable separatrices $S_\mu^+$ and $S_\mu^-$, respectively. Consider two closed connected arcs $l^{\pm}\subset S^{\pm}_{\mu_0}$, having both an endpoint at $p_{\mu_0}$. In the case of a homoclinic connection ($S_{\mu_0}^+=S_{\mu_0}^-$), we require additionally that $l^+\cap l^-=\{p_{\mu_0}\}$. Then, there exists a neighborhood $V$ of $(l^+\cup l^-)\times{\mu_0}$ in $\mathbb{R}^2\times\mathbb{R}^N$ and a $\mathscr C^\infty$ diffeomorphism $\Phi:V\to\Phi(V)\subset\mathbb{R}^2\times\mathbb{R}^N$ with $\Phi(x,y,\mu)=(\phi_\mu(x,y),\mu)$ such that $(\phi_\mu)_{*}(X_\mu)=xP(x,y;\mu)\partial_x+yQ(x,y;\mu)\partial_y$, with $P,Q\in \mathscr C^\infty(\Phi(V))$.
\end{lem}

\section*{Acknowledgments}

This work is financially supported by the Spanish Ministry of Science, Innovation and Universities, through grants PID2021-125625NB-I00 and PID2020-118281GB-C33 and by the Agency for Management of University and Research Grants of Catalonia through grants 2021SGR01015 and 2021SGR00113. This work is also supported by the Spanish State Research Agency, through the Severo Ochoa and Mar\'ia de Maeztu Program for Centers and Units of Excellence in R\&D (CEX2020-001084-M). The second author is supported by S\~ao Paulo Research Foundation (FAPESP) grants 21/14450-4 and 19/13040-7.

\bibliographystyle{siam}
\bibliography{Referencias_BEPE.bib}

\end{document}